\documentclass[a4paper,11pt]{article}
\usepackage{a4wide}
\usepackage{theorem}
\usepackage{amsmath}
\usepackage{array}
\usepackage{amssymb}
\usepackage{amsfonts}
\usepackage[french,english]{babel}
\usepackage{epsf}
\usepackage{epsfig}

\newtheorem{theo}{\indent Theorem\newline}[section]

{\theorembodyfont{\rmfamily}%
\newtheorem{rem}[theo]{\noindent Remark}}
{\theorembodyfont{\rmfamily}%
 \theoremstyle{break}%
}
\newtheorem{prop}[theo]{\indent Proposition\newline}
\newtheorem{lemma}[theo]{\indent Lemma\newline}
\newtheorem{cor}[theo]{\indent Corollary\newline}

 \def\N{{\Bbb{N}}}
\def\Z{{\Bbb{Z}}}

\def\R{{\Bbb{R}}}
\def\C{{\Bbb{C}}}

\usepackage{eufrak}
\newcommand{\goth}[1]{\EuFrak{#1}}
\newcommand{\s}{\mathop{\goth s}\nolimits}
\newcommand{\p}{\mathop{\goth p}\nolimits}

\newcommand{\inv}{\mathop{\rm inv}\nolimits}

\newcommand{\coker}{\mathop{\rm coker}\nolimits}

\newlength{\indentation}%
\makeatletter
\newcommand\@makefntextsans[1]{%
    \parindent 0em%
    \noindent%
    \hb@xt@0em{\hss}%
    #1}
\def\footnotetextsans{%
     \@ifnextchar [\@xfootnotenextsans%
       {\@footnotetextsans}}
\def\@xfootnotenextsans[#1]{%
  \begingroup%
     \csname c@\@mpfn\endcsname #1\relax%
  \endgroup%
  \@footnotetextsans}
\long\def\@footnotetextsans#1{\insert\footins{%
    \reset@font\footnotesize%
    \interlinepenalty\interfootnotelinepenalty%
    \splittopskip\footnotesep%
    \splitmaxdepth \dp\strutbox \floatingpenalty \@MM%
    \hsize\columnwidth \@parboxrestore%
    \color@begingroup%
      \@makefntextsans{%
        \rule\z@\footnotesep\ignorespaces#1\@finalstrut\strutbox}
    \color@endgroup}}
\makeatother

\begin{document}

\cleardoublepage
\title{Spinor states of real rational curves in real algebraic convex $3$-manifolds
and enumerative invariants}
\author{Jean-Yves Welschinger\footnote{Member of the european network RAAG CT-2001-00271}}
\maketitle

\makeatletter\renewcommand{\@makefnmark}{}\makeatother
\footnotetextsans{Keywords :  Convex manifold, real algebraic manifold, 
stable map, enumerative geometry.}
\footnotetextsans{AMS Classification : 14N35, 14P25 .}

{\bf Abstract :}

Let $X$ be a real algebraic convex $3$-manifold whose real part is equipped
with a $Pin^-$ structure. We show that every irreducible real rational curve
with non-empty real part has a canonical spinor state belonging to $\{\pm 1\}$. 
The main result is then that the algebraic count of the number of 
real irreducible rational curves in a given numerical equivalence class passing
through the appropriate number of points does not depend on the choice of the real
configuration
of points, provided that these curves are counted with respect to their spinor states.
These invariants provide lower bounds for the total number of such real rational 
curves independantly of the choice of the real configuration of points.

\section*{Introduction}

A smooth complex algebraic projective manifold $X$ is said to be {\it convex} when
the vanishing $H^1 (\C P^1 ; u^* TX)=0$ occurs for every morphism $u : \C P^1 \to X$.
Main examples are homogeneous spaces. These manifolds 
provide a suitable framework in order to define genus $0$ algebraic Gromov-Witten 
invariants
since the space of morphisms from $\C P^1$ to $X$ in a given homology class 
$d \in H_2 (X ; \Z)$ turns out to be a smooth manifold of the expected 
dimension $c_1 (X)d + 3$ (see \cite{Kont}, \cite{FP}). In particular, these Gromov-Witten 
invariants are enumerative. Let $X$ be a convex manifold of dimension $3$ and
 $d \in H_2 (X ; \Z)$ be such that $c_1 (X)d$ is even. Let $k_d$ be half of this even
integer. Then, through a generic configuration $\underline{x} = (x_1 , \dots , x_{k_d})$
of $k_d$ distinct points of $X$ passes only finitely many connected rational curves
in the homology class $d$. These curves are all irreducible and immersed and their
number $N_d$ does not depend on the choice of $\underline{x}$. It is a Gromov-Witten
invariant of the manifold $X$, often denoted by $GW(X,d,pt, \dots, pt)$. From now on,
assume that $X$ is {\it real}, that is equipped with an antiholomorphic involution $c_X$.
The fixed point set of $c_X$ is denoted by $\R X$ and called the {\it real part} of $X$.
Assume also that $(c_X)_* d = -d$ and that the configuration $\underline{x}$ is {\it real}, 
that is satisfies $\{ c_X (x_1) , \dots , c_X (x_{k_d}) \} = \{ x_1 , \dots , x_{k_d} \}$. 
Then, the analogous result is no more true, 
the number $R_d (\underline{x})$ of irreducible rational curves in the homology class
$d$ passing through $\underline{x}$ that are real depends in general on the choice of 
$\underline{x}$. The parity of this integer is however invariant,
it is the one of $N_d$. The main purpose of this paper is to refine this $mod(2)$
invariant into an integer valued one, see Theorem \ref{maintheorem}.

Fix a $Pin^-_3$ structure $\p$ on $\R X$ (see \S \ref{subsectspinstruct}) and assume that 
$\underline{x} \cap \R X \neq \emptyset$. As soon as 
the choice of $\underline{x}$ is generic enough, this implies that for every
curve $A \in {\cal R}_d (\underline{x})$, where ${\cal R}_d (\underline{x})$ denotes
the set of real irreducible rational curves in the homology class $d$ passing through 
$\underline{x}$, the real part $\R A$ is nonempty. Also, as soon as 
$\underline{x}$ is generic enough, the normal bundle $N_A$ of $A$ in $X$
is real and {\it balanced}, that is the direct sum of two isomorphic holomorphic
line bundles. Choosing
a real subline bundle of $N_A$ of maximal degree $k_d - 1$, we can equip the knot
$\R A \subset \R X$ with a canonical ribbon structure, see \S \ref{subsectspinorient} for
the detailed construction. This structure does depend on the way the complex curve $A$
is immersed in the complex manifold $X$ and it is inherited from this complex immersion.
Associating a framing to this ribbon knot, we construct a loop in the 
$O_3 (\R)$-principal bundle $R_X$ of orthonormal frames of $\R X$. The
{\it spinor state} $sp(A) \in \{ \pm 1 \}$ of $A$ is then defined to be the
obstruction to lift this loop into a loop of the $Pin^-_3$-principal bundle $P_X$
given by the pin structure $\p$. Denote by $(\R X)_1, \dots , (\R X)_n$ the connected
components of $\R X$ and, for $i \in \{ 1 , \dots , n \}$, by $r_i = \# (\underline{x}
\cap (\R X)_i)$. Then define 
$$\chi^{d,\p}_r (\underline{x}) = \sum_{A \in {\cal R}_d (\underline{x})} sp(A),$$
where $r = (r_1 , \dots , r_n)$.
\begin{theo}
\label{theoprincintro}
1) As soon as $c_1(X)d \neq 4$, the integer $\chi^{d,\p}_r (\underline{x})$ is
independant of the choice of the real configuration $\underline{x} \in X^{k_d}$, it only
depends on $d, \p$ and the $n$-tuple $r$. The same holds when $c_1(X)d = 4$ and
$X$ does not have any non-immersed rational curve in the class $d$.

2)The similar construction and result hold for the blown up projective space
$Y = \C P^3 \# \overline{\C P}^3$ when the homology class $d$ satisfies 
$Exc. d \geq 0$, where $Exc$ is the exceptional divisor of $Y$.
\end{theo}
Note that the blown up projective space is not convex, see \S \ref{sectfurth} for the
second part of this theorem \ref{theoprincintro}. See also Theorem \ref{maintheorem}
for a slightly more general statement in the case $c_1(X)d = 4$.

This integer $\chi^{d,\p}_r (\underline{x})$ is denoted by $\chi^{d,\p}_r$ and we put
$\chi^{d,\p}_r = 0$ as soon as it is not well defined. This allows to define the
polynomial $\chi^{d,\p} (T) = \sum_{r \in \N^n} \chi^{d,\p}_r T^r \in \Z[T]$, where
$T^r = T_1^{r_1} \dots T_n^{r_n}$. This polynomial is of the same parity as the integer
$k_d$ and each of its monomials actually only depends on one indeterminate. Theorem
\ref{theoprincintro} means that the function
$\chi^{\p} : d \in H_2 (X ; \Z) \mapsto \chi^{d,\p} (T) \in \Z[T]$ is an invariant
associated to the isomorphism class of the real algebraic convex manifold $(X , c_X)$.
As an application, this invariant provides the following lower bounds in real enumerative
 geometry.

\begin{cor}
\label{corintro}
Under the assumptions of Theorem \ref{theoprincintro}, we have 
$$|\chi^{d,\p}_r| \leq R_d (\underline{x}) \leq N_d,$$
independantly of the choice of $\underline{x} \in X^{k_d}$. 
\end{cor}

Note that a similar invariant and similar lower bounds have already been obtained in 
\cite{Wels}, \cite{Wels2}  using real rational curves in real symplectic $4$-manifolds.
The question was then raised whether there exists such invariants in higher dimensions.
The results presented here thus provide a partial answer to this question.

Let us fix now $(X , c_X) = (\C P^3 , conj)$ and a spin structure $\s$ on $\R P^3$. Note that
for orientable components of $\R X$, a spin structure is more convenient for us than a $Pin^-$
structure, see Remark \ref{centrem}. The $n$-tuple $r$ is then reduced to a
single even integer between $2$ and $2d = k_d$. Denote by 
$Y = \C P^3 \# \overline{\C P}^3$
the blown up of $X$ at a real point $x_0$. Denote by $l$ the homology class of a line
in the exceptional divisor $Exc$ of $Y$ and by $f$ the class of the strict transform of 
a line in $\C P^3$ passing through $x_0$. Set $d_Y = d(f+l) - 2l$ and choose a spin 
structure $\overline{\s}$ on $\overline{\R P}^3$. It induces a spin structure
$\s \# \overline{\s}$ on $\R Y = \R P^3 \# \overline{\R P}^3$, see \S \ref{subsectfurth}.
From the second part of Theorem \ref{theoprincintro}, the integer 
$\chi_r^{d_Y, \s \# \overline{\s}}$ is well defined for $r$ even between $2$ and 
$k_d - 2$ and
it is an invariant associated to $Y$ together with its real structure. The following 
theorem provides relations between the coefficients of the polynomial $\chi^{d,\s} (T)$
for the projective space, see Theorem \ref{theorel}.

\begin{theo}
\label{introrel}
Let $r$ be an even integer between $2$ and $k_d - 2$, then
$$\chi^{d,\s}_{r+2} = \chi^{d,\s}_r - 2\chi_r^{d_Y, \s \# \overline{\s}}.$$
\end{theo}

This paper is divided in four paragraphs. In the first one, we introduce Kontsevich's
space of stable maps and its real structures, as well as preliminaries on the evaluation
map and Gromov-Witten invariants. In the second one, we introduce the definition of
spinor states and state the main results of this paper. In the third one, we give
a proof of these results. Finally, the last paragraph is devoted to a further study
of the polynomial $\chi^{d,\s} (T)$ for the projective space, proving
Theorem \ref{introrel}.\\

{\bf Acknowledgements :}

I am grateful to V. Kharlamov and O. Viro for fruitful discussions on spin and pin structures.
Part of this work has been done during my stay at the Universit\'e Louis Pasteur of
Strasbourg in July $2003$. I would like to acknowledge the warm hospitality of people working 
there.

\tableofcontents

\section{Evaluation map and Gromov-Witten invariants}

Let $(X,c_X)$ be a real algebraic convex $3$-manifold and $d \in H_2 (X ; \Z)$ be such
that $(c_X)_* d = -d$.

\subsection{Moduli space of genus $0$ real stable maps with $k$ marked points}
\label{subsectmoduli}

Let ${\cal M}or_d (X)$ be the set of morphisms $u$ from $\C P^1$ to $X$ {\it in the 
class $d$},
that is satisfying $u_* [\C P^1] = d$. It is a smooth quasi-projective manifold of
pure
dimension $c_1 (X)d +3$. Let $k \in \N^*$ and $Diag_k = \{ (z_1 , \dots , z_k) \in 
(\C P^1)^k \,
| \, \exists \, i \neq j \, , \,  z_i = z_j \}$. The manifold ${\cal M}or_d (X)
\times ((\C P^1)^k \setminus Diag_k)$ is equipped with an action of the Moebius group
defined by :
$$
\begin{array}{rcl}
{\cal M}oeb \times {\cal M}or_d (X)
\times ((\C P^1)^k \setminus Diag_k) & \to & {\cal M}oeb \times {\cal M}or_d (X)
\times ((\C P^1)^k \setminus Diag_k) \\
(\phi , (u,z_1 , \dots , z_k)) & \mapsto & \left\{ 
\begin{array}{l}
(u \circ \phi^{-1} , \phi (z_1) , \dots , \phi (z_k)) \text{ if } \phi \in 
Aut (\C P^1),\\
(c_X \circ u \circ \phi^{-1} , \phi (z_1) , \dots , \phi (z_k)) \text{ otherwise.}
\end{array}
\right.
\end{array}$$
It is also equipped with an action of the symmetric group ${\cal S}_k$ defined by
$(\sigma , (u,z_1 , \dots , z_k)) \mapsto (u,z_{\sigma(1)} , \dots , z_{\sigma(k)})$
for $(\sigma , (u,z_1 , \dots , z_k)) \in {\cal S}_k \times {\cal M}or_d (X)
\times ((\C P^1)^k \setminus Diag_k)$. Denote by ${\cal M}^d_{0,k} (X)$ the quotient
of ${\cal M}or_d (X)
\times ((\C P^1)^k \setminus Diag_k)$ by $Aut (\C P^1)$. This moduli space is a 
quasi-projective
manifold of pure dimension $c_1 (X)d + k$ and a complex orbifold. It is
equipped with an action of ${\cal S}_k$ and with an antiholomorphic involution 
$c_{\cal M}$ induced by the action of
${\cal M}oeb / Aut (\C P^1) \cong \Z / 2\Z$.
Denote by ${\cal M}^d_{0,k} (X)^*$ the complement of the set of maps $u$ which
can be written $u' \circ g$ where $u' : \C P^1 \to X$ and
$g : \C P^1 \to \C P^1 $ is a non-trivial ramified covering. It is included
in the smooth locus of ${\cal M}^d_{0,k} (X)$. Denote by
$U^d_{0,k} (X) \to {\cal M}^d_{0,k} (X)^*$ the universal curve, which is obtained as
the quotient of ${\cal M}or_d (X)
\times ((\C P^1)^k \setminus Diag_k) \times \C P^1$ by the action of $Aut (\C P^1)$ and
then restricted over ${\cal M}^d_{0,k} (X)^*$. This universal curve is also equipped 
with an
antiholomorphic involution $c_{\cal U}$ and with an action of ${\cal S}_k$ which
lift $c_{\cal M}$ and the action of ${\cal S}_k$ on ${\cal M}^d_{0,k} (X)^*$ respectively.
Finally, denote by $\overline{{\cal M}}^d_{0,k} (X)$ the moduli space of genus $0$
stable maps with $k$ marked points in the homology class $d$. This is the compactification
due to M. Kontsevich of the space ${\cal M}^d_{0,k} (X)$, see \cite{Kont}, \cite{FP}. Denote by 
$\overline{{\cal M}}^d_{0,k} (X)^*$ the space of stable maps $(u, C ,z_1 , \dots , z_k)$
for which $u$ restricted to any irreducible component of $C$ does not
admit a factorization of the form $u' \circ g$ where $u' : \C P^1 \to X$ and
$g : \C P^1 \to \C P^1 $ is a non-trivial ramified covering. Note that this in
particular avoids maps
$(u, C ,z_1 , \dots , z_k)$ mapping a component of $C$ to a constant, which does 
not matter for us since these maps will play no r\^ole in this paper. We recall the following
theorem.

\begin{theo}
\label{theoKont}
1) The manifold $\overline{{\cal M}}^d_{0,k} (X)$ is projective normal of pure
dimension $c_1 (X)d + k$. It is a complex orbifold containing 
$\overline{{\cal M}}^d_{0,k} (X)^*$ in its smooth locus, which is equipped with an
evaluation morphism $ev^d_k : \overline{{\cal M}}^d_{0,k} (X) \to X^k$. Finally,
the complement $\overline{{\cal M}}^d_{0,k} (X)^* \setminus {\cal M}^d_{0,k} (X)^*$
is a divisor with normal crossings.

2)The curve $U^d_{0,k} (X) \to {\cal M}^d_{0,k} (X)^*$ extends to a universal curve
$\overline{U}^d_{0,k} (X) \to \overline{{\cal M}}^d_{0,k} (X)^*$.

3)The real structures $c_{\cal M}$ and $c_{\cal U}$ extend to real structures
$c_{\overline{{\cal M}}}$ and $c_{\overline{{\cal U}}}$ on 
$\overline{{\cal M}}^d_{0,k} (X)$
and $\overline{U}^d_{0,k} (X)$ respectively. The same holds for the actions of 
${\cal S}_k$ on 
these spaces. $\square$
\end{theo}
This theorem is proved in \cite{FP}, Theorems $2$ and $3$. The third part of this theorem
follows from the fact that the construction presented in \cite{FP} can be carried out
over the reals, see \cite{Kwon}. Note that the evaluation morphism restricted to
${\cal M}^d_{0,k} (X)$ is just the map induced by
$(u, z_1 , \dots , z_k) \in {\cal M}or_d (X)
\times ((\C P^1)^k \setminus Diag_k) \mapsto (u(z_1) , \dots , u(z_k)) \in X^k$.

For every element $\tau \in {\cal S}_k$, denote by $\Phi_{\tau}$ the associated
automorphism of $\overline{{\cal M}}^d_{0,k} (X)$. This automorphism commutes with
$c_{\overline{{\cal M}}}$. When $\tau^2 = id$, put $c_{\overline{{\cal M}}, \tau}
=\Phi_{\tau} \circ c_{\overline{{\cal M}}}$ the associated real structure. Note that
if $\tau' = \sigma \circ \tau \circ \sigma^{-1}$, with $\tau', \sigma \in {\cal S}_k$,
then $c_{\overline{{\cal M}}, \tau'} = \Phi_{\sigma} \circ c_{\overline{{\cal M}}, \tau}
\circ \Phi_{\sigma}^{-1}$. Thus, the isomorphism class of 
$(\overline{{\cal M}}^d_{0,k} (X) , c_{\overline{{\cal M}}, \tau})$ only depends on 
the conjugacy class of $\tau$ in ${\cal S}_k$. Similarly, the group ${\cal S}_k$ acts
on $X^k$ by permutation of the factors. Denote by $c_{id}$ the real structure
$c_X \times \dots \times c_X$ on $X^k$. The evaluation morphism $ev^d_k$ is
obviously equivariant for the actions of ${\cal S}_k$ and $\Z / 2\Z$. For every element
$\tau \in {\cal S}_k$ such that $\tau^2 = id$, we set $c_{\tau}$ to be the real
structure $\Psi_{\tau} \circ c_{id}$ on $X^k$, where $\Psi_{\tau}$ is the automorphism
of $X^k$ associated to $\tau$. Then, the isomorphism class of $(X^k , c_{\tau})$ only
depends on the conjugacy class of $\tau$ in ${\cal S}_k$ and
$ev^d_k : (\overline{{\cal M}}^d_{0,k} (X) , c_{\overline{{\cal M}}, \tau}) \to 
(X^k , c_{\tau})$ is a real morphism between real algebraic varieties.

\subsection{Evaluation map and balanced curves}

\subsubsection{Balanced rational curves}
\label{subsectbalanced}

Let $(u, z_1 , \dots , z_k) \in {\cal M}or_d (X)
\times ((\C P^1)^k \setminus Diag_k)$. The differential of $u$ induces a morphism 
of sheaves $0 \to {\cal O}_{\C P^1} (T \C P^1) \to {\cal O}_{\C P^1} (u^* TX)$.
Denote by ${\cal N}_u$ the quotient sheaf. This quotient sheaf admits a decomposition
${\cal O}_{\C P^1} (N_u) \oplus {\cal N}_u^{sing}$, where $N_u$ is the normal bundle
of $u$ and ${\cal N}_u^{sing}$ is the skyscraper sheaf having the critical points
$p_i \in \C P^1$ of $du$ as support, and $\C^{n_i}$ as fibers, where $n_i$ is the
vanishing order of $du$ at $p_i$. The bundle $N_u$ is holomorphic of rank two over
$\C P^1$. From a theorem of Grothendieck, it splits as
the direct sum of two line bundles ${\cal O}_{\C P^1} (a) \oplus {\cal O}_{\C P^1} (b)$, 
see \S $2.1$ of \cite{OSS} for example.
This normal bundle is said to be {\it balanced} when $a=b$. The rational curve
$(u, z_1 , \dots , z_k)$ is said to be {\it balanced} when $u$ is an immersion and
$N_u$ is balanced. Note that when $u$ is an immersion, the adjunction formula imposes
$a+b = c_1 (X)d -2$. Thus a necessary condition for the curve to be balanced is that
$c_1 (X)d$ is even.\\

{\bf Examples :}

1) The rational curve $u : (t_0 , t_1) \in \C P^1 \mapsto (t_0^d , t_0^{d-1} t_1 ,
t_0 t_1^{d-1} , t_1^d) \in \C P^3$ is balanced as soon as $d \in \N^*$ is different
from  $2$.

2) If $A \subset X^3$ is balanced and $Y$ is the blown up of $X^3$ at a point
$x_0 \in A$, then the strict transform of $A$ in $Y$ is balanced, see Lemma 
\ref{lemmablowup}.

\subsubsection{Dominance of the evaluation map}

From now on, we assume that $c_1 (X)d$ is even and we set $k_d = \frac{1}{2} c_1 (X)d$.
Note that this condition is automatically fulfilled when $\R X$ is orientable and
$d$ is realized by a real $2$-cycle $A$, since then $c_1 (X)[A] = w_1 (\R X) [\R A] = 0
\mod(2)$.

Let $(u, \C P^1 , z_1 , \dots , z_{k_d}) \in {\cal M}^d_{0,k_d} (X)^*$ and set
$N_{u , -z} = N_u \otimes {\cal O}_{\C P^1} (-z)$, where $z = ( z_1, \dots , z_{k_d} )
\in (\C P^1)^{k_d}$. Then, set ${\cal N}_{u , -z} = {\cal O}_{\C P^1} (N_{u , -z}) 
\oplus {\cal N}_u^{sing}$.

\begin{lemma}
\label{lemmacoker}
Let $(u, \C P^1 , z_1 , \dots , z_{k_d}) \in {\cal M}^d_{0,k_d} (X)^*$, the following
isomorphisms occur :
$$\ker (d|_{(u, \C P^1 , z)} ev^d_{k_d}) \cong H^0 (\C P^1 ; {\cal N}_{u , -z}) \cong
H^0 (\C P^1 ; N_{u , -z}) \oplus H^0 (\C P^1 ; {\cal N}_u^{sing}),$$
$$\text{and } \coker (d|_{(u, \C P^1 , z)} ev^d_{k_d}) \cong H^1 (\C P^1 ; N_{u , -z}).$$
\end{lemma}
In particular, balanced curves are regular points of the evaluation
map $ev^d_{k_d}$.\\

{\bf Proof :}

Denote by $eval^d_{k_d}$ the map $(u, z_1 , \dots , z_{k_d}) \in {\cal M}or_d (X)
\times ((\C P^1)^{k_d} \setminus Diag_{k_d}) \mapsto (u(z_1) , \dots , u(z_{k_d})) 
\in X^{k_d}$.
The differential of this map at the point $(u,z)$ is given by
$(v, \stackrel{.}{z}) \in H^0 (\C P^1 ; u^* TX) \times T_z (\C P^1)^{k_d} \mapsto
v(z) + d|_z u (\stackrel{.}{z}) \in T_{u(z)} X^{k_d}$. Denote by ${\cal N}_{u|_z}$
the skyscraper sheaf with support $\{ z_1 , \dots , z_{k_d} \}$ and fiber 
$T_{u(z_i)} X / d|_{z_i} u (T \C P^1)$ over $z_i$.
The cokernel of $\stackrel{.}{z} \in T_z (\C P^1)^{k_d} \mapsto d|_z u (\stackrel{.}{z})
\in T_{u(z)} X^{k_d}$ is isomorphic to $H^0 (\C P^1 ; {\cal N}_{u|_z})$. Thus, the 
cokernel of $d|_{(u, z)} eval^d_{k_d}$ is identified with the cokernel of the 
composition :
$$\begin{array}{ccccc}
H^0 (\C P^1 ; u^* TX) & \to & H^0 (\C P^1 ; {\cal N}_u) & \to & H^0 (\C P^1 ; 
{\cal N}_{u|_z})
\quad (*)\\
v & \mapsto & [v] & \mapsto & [v(z)].
\end{array}$$
From the exact sequence of sheaves $0 \to T \C P^1 \to u^*TX \to {\cal N}_u \to 0$
we deduce that the first morphism of $(*)$ is surjective since $H^1 (\C P^1 ; T \C P^1)
=0$ and that $H^1 (\C P^1 ; N_u) = 0$ since $X$ is convex. From the short exact 
sequence of sheaves $0 \to {\cal N}_{u , -z} \to {\cal N}_u \to {\cal N}_{u|_z} \to 0$,
we deduce the long exact sequence 
$$\dots \to H^0 (\C P^1 ; {\cal N}_u) \to H^0 (\C P^1 ;{\cal N}_{u|_z} ) \to
H^1 (\C P^1 ;{\cal N}_{u , -z} ) \to H^1 (\C P^1 ; {\cal N}_u) \to \dots$$
Since $ H^1 (\C P^1 ; {\cal N}_u) =0$, the cokernel of the second morphism
of $(*)$ is isomorphic to $H^1 (\C P^1 ;{\cal N}_{u , -z} )$. Thus, the cokernel
of $d|_{(u, z)} eval^d_{k_d}$ and hence the cokernel of $d|_{(u, \C P^1 , z)} ev^d_{k_d}$
is isomorphic to $H^1 (\C P^1 ;{\cal N}_{u , -z} )$.

In the same way, the kernel of the map $d|_{(u, \C P^1 , z)} ev^d_{k_d}$ is isomorphic to
the quotient of the kernel of the composition $(*)$ with the image of the
morphism $du : H^0 (\C P^1 ; T \C P^1) \to H^0 (\C P^1 ; u^*TX )$. Indeed, the latter
coincide with the infinitesimal action of $Aut (\C P^1)$ on ${\cal M}or_d (X)$.
Since the quotient of $H^0 (\C P^1 ; u^*TX )$ by $du( H^0 (\C P^1 ; T \C P^1))$
is isomorphic to $ H^0 (\C P^1 ; {\cal N}_u)$, we deduce that the kernel of
$d|_{(u, \C P^1 , z)} ev^d_{k_d}$ is isomorphic to the kernel of the second
morphism of $(*)$, that is to $H^0 (\C P^1 ; {\cal N}_{u , -z})$. $\square$

\begin{prop}
\label{propcoker}
Let $(u , C ,  z) \in \overline{{\cal M}}^d_{0,k_d} (X)^*$ be such 
that $C$ has two irreducible components $C_1$ and $C_2$ for which $u(C_1)$, $u(C_2)$ are
balanced and meet in a single ordinary double point away from 
$u(z)$. Assume moreover that for $i \in \{ 1,2\}$,
$\# (C_i \cap z) = \frac{1}{2} c_1(X)d_i$ where $d_i = u_* [C_i] \in H_2 (X ; \Z)$. Then,
$(u, C , z)$ is a regular point of $ev^d_{k_d}$.
\end{prop}

{\bf Proof :}

For $i \in \{ 1,2\}$, set $k_{d_i} = \frac{1}{2} c_1(X)d_i$, $u_i = u|_{C_i}$ and 
$z^i = C_i \cap z$.
Denote by ${\cal K} = \overline{{\cal M}}^d_{0,k_d} (X)^* \setminus 
{\cal M}^d_{0,k_d} (X)^*$ the divisor of stable maps having reducible domain. Then
$(u, C , z)$ is a smooth point of ${\cal K}$. Indeed, for $i \in \{ 1,2\}$, denote
by $e_i : \overline{{\cal M}}^{d_i}_{0,k_{d_i} \cup \{ \bullet \}} (X) \to X$ the
evaluation map associated to the additional marked point $\bullet$. Since
$u(C_i)$ is balanced, $d|_{(u_i, C_i , z^i \cup (C_1 \cap C_2))} e_i$ is surjective. 
Thus, $\big( (u_1, C_1 , z^1 \cup (C_1 \cap C_2)) , (u_2, C_2 , z^2 \cup 
(C_1 \cap C_2)) \big)$ is a smooth point of $\widetilde{\cal K} = (e_1 \times e_2)^{-1} (Diag_X)$,
where $Diag_X \subset X \times X$ is the diagonal. From Lemma $12 (i)$ of \cite{FP},
the map $\widetilde{\cal K} \to {\cal K}$ is an isomorphism, hence $(u, C , z)$ is a smooth 
point of ${\cal K}$.

Now, the restriction of $ev^d_{k_d}$ to ${\cal K}$ in the neighborhood of $(u, C , z)$ writes
as the composition of the morphism of restriction $(u, C_1 \cup C_2 , z) \in {\cal K}
\mapsto ((u_1, C_1 , z^1) , (u_2, C_2 , z^2 )) \in {\cal M}^{d_1}_{0 , k_{d_1}} (X)
\times {\cal M}^{d_2}_{0 , k_{d_2}} (X)$ and the morphism of evaluation
$ev^{d_1}_{k_{d_1}} \times ev^{d_2}_{k_{d_2}} : {\cal M}^{d_1}_{0 , k_{d_1}} (X)
\times {\cal M}^{d_2}_{0 , k_{d_2}} (X) \to X^{k_d}$. Since $u(C_1)$ and $u(C_2)$ meet
in only one ordinary double point, the first morphism is injective in a neighborhhod of
$(u, C , z)$. Since $u(C_1)$ and $u(C_2)$ are balanced, from Lemma \ref{lemmacoker},
the second morphism is an isomorphism in the neighborhood of $((u_1, C_1 , z^1) , 
(u_2, C_2 , z^2 ))$. Hence, to prove proposition \ref{propcoker}, it suffices to
prove that $ev^d_{k_d}$ is injective when restricted to a transversal of ${\cal K}$ at
the point $(u, C , z)$. Since $ev^d_{k_d} ({\cal K})$ is smooth in a neighborhood of $u(z)$,
a transversal to this divisor can be chosen in one factor of $X^{k_d}$, that is of the
form $\gamma : t \in \Delta \mapsto (x_1 (t), \dots , x_{k_d} (t)) \in X_{k_d}$ where only
one point $x_i (t)$, say $x_{k_d} (t)$, indeed depends on $t$. Since such a path is transversal
to $ev^d_{k_d}$, its inverse image is a smooth curve $B$ in 
$ \overline{{\cal M}}^d_{0,k_d} (X)^*$ transversal to ${\cal K}$ at the point $(u, C , z)$.
Denote by $U \to B$ the restriction of $\overline{U}^d_{0 , k_d} (X) \to
\overline{{\cal M}}^d_{0,k_d} (X)^*$ over $B$. This universal curve has one unique reducible
fiber over $(u, C , z)$ (choosing a smaller $B$ if necessary) and $k_d$ tautological
sections $s_1 , \dots , s_{k_d} : B \to U$. The evaluation morphism $ev_U : U \to X$ contracts
the curves $Im (s_1) , \dots ,  Im (s_{k_d - 1})$ and we have to prove that it is injective
when restricted to $Im (s_{k_d})$. This follows from the following Lemma \ref{lemmaevU}.
$\square$

\begin{lemma}
\label{lemmaevU}
The differential of the evaluation morphism $ev_U : U \to X$ defined above is injective at
every point of the central fiber $C \setminus \{ z_1 , \dots , z_{k_d - 1} \}$.
\end{lemma}

{\bf Proof :}

Without loss of generality, we can assume that $z^1 = \{ z_1 , \dots , z_{k_{d_1}} \}$ and
$z^2 = \{ z_{k_{d_1} + 1} , \dots , z_{k_{d}} \}$. Since $u|_{C_1}$ and $u|_{C_2}$ are
immersions, the rank of $d ev_U$ is at least one at each point of $C$. Now, the normal bundle of
$C_1$ in $U$ is isomorphic to ${\cal O}_{C_1} (-1)$. The morphism $d ev_U$ thus induces
a morphism of sheaves $0 \to {\cal O}_{C_1} (-1) \to {\cal O}_{C_1} (N_{u_1})$ where
$u_1 = u|_{C_1}$. Since this morphism vanishes at the points $z_1 , \dots , z_{k_{d_1}}$, its
image is a subline bundle of $N_{u_1}$ of degree at least $k_{d_1} - 1$. From the isomorphism
$N_{u_1} \cong {\cal O}_{C_1} (k_{d_1} - 1) \oplus {\cal O}_{C_1} (k_{d_1} - 1)$, we deduce
that this subline bundle has degree exactly $k_{d_1} - 1$. This implies that the morphism
of sheaves vanishes nowhere else than at the points $z_1 , \dots , z_{k_{d_1}}$ and thus
$d ev_U$ is of rank two at each point of $C_1 \setminus z^1$. Similarly, the normal bundle
of $C_2$ in $U$ is isomorphic to ${\cal O}_{C_2} (-1)$. The morphism $d ev_U$ thus induces
a morphism of sheaves $0 \to {\cal O}_{C_2} (-1) \to {\cal O}_{C_2} (N_{u_2})$, where
$u_2 = u|_{C_2}$, which vanishes at the points $z_{k_{d_1} + 1} , \dots , z_{k_{d} - 1}$.
The image of this morphism is thus a subline bundle $L$ of $N_{u_2}$ of degree at least 
$k_{d_2} - 2$. Denote by $L_1$ the unique subline bundle of $N_{u_2} \cong {\cal O}_{C_2} 
(k_{d_2} - 1) \oplus {\cal O}_{C_2} (k_{d_2} - 1)$ which is of degree $k_{d_2} - 1$ and
contains the tangent line of $u (C_1)$ at $u (C_1) \cap u (C_2)$. Denote by $M_1$ the
quotient $N_{u_2} / L_1$, it is a line bundle of degree $k_{d_2} - 1$. The bundle $L_1$
corresponds to infinitesimal deformations of the curve $u (C_2)$ into a curve passing
through $u (z_{k_{d_1} + 1}) , \dots , u (z_{k_{d} - 1})$ and $u (C_1)$. Since by construction
$B$ is transversal to the divisor ${\cal K}$ of reducible curves of 
$\overline{{\cal M}}^d_{0,k_d} (X)^*$, one has $L \neq L_1$. However, $L$ also contains
the tangent line of $u (C_1)$ at $u (C_1) \cap u (C_2)$. It follows that the morphism of
sheaves $0 \to {\cal O}_{C_2} (L) \to {\cal O}_{C_2} (M_1)$ induced by the projection
$N_{u_2} \to M_1$ vanishes at the point $u (C_1) \cap u (C_2)$. This implies that
$\deg (L) \leq \deg (M_1) - 1 = k_{d_2} - 2$. Hence, the line bundle $L$ is of degree
exactly $k_{d_2} - 2$ and $d ev_U$ is of rank two at each point of $C_2 \setminus \{
z_{k_{d_1} + 1} , \dots , z_{k_{d} - 1} \}$. The result follows. $\square$

\begin{prop}
\label{propcoker2}
Let $(u , C , z) \in \overline{{\cal M}}^d_{0,k_d} (X)^*$  be such 
that $C$ has two irreducible components $C_1$ and $C_2$ for which $u(C_1)$, $u(C_2)$ are
immersed and meet in a single ordinary double point away from 
$u(z)$. For $i \in \{ 1,2\}$, denote by 
$z^i = C_i \cap z$, $d_i = u_* [C_i] \in H_2 (X ; \Z)$,
$k_{d_1} = E( \frac{1}{2} c_1(X)d_1 )$ and $k_{d_2} = E( \frac{1}{2} c_1(X)d_2 ) + 1$ where
$E()$ denotes the integer part. Assume that $N_{u|_{C_1}} \cong {\cal O}_{C_1} 
(k_{d_1} - 1) \oplus {\cal O}_{C_1} (k_{d_1})$, $N_{u|_{C_2}} \cong {\cal O}_{C_2} 
(k_{d_2} - 2) \oplus {\cal O}_{C_2} (k_{d_2} - 1)$ and that the tangent line
to $u(C_1)$ (resp. to $u(C_2)$) at the point $u(C_1) \cap u(C_2)$ is not mapped to the unique
subline bundle of degree $k_{d_2} - 1$ (resp. $k_{d_1}$) of $N_{u|_{C_2}}$ 
(resp. $N_{u|_{C_1}}$). Then,
$(u, C , z)$ is a regular point of $ev^d_{k_d}$.
\end{prop}

{\bf Proof :}

As in the proof of proposition \ref{propcoker}, the divisor ${\cal K} = 
\overline{{\cal M}}^d_{0,k_d} (X)^* \setminus {\cal M}^d_{0,k_d} (X)^*$ is smooth in a 
neighborhood of $(u, C , z)$. Moreover, $ev^d_{k_d}$ restricted to ${\cal K}$ is injective
in the neighborhood of $(u, C , z)$.
Indeed, this restriction writes
as the composition of the morphism $(u, C_1 \cup C_2 , z) \in {\cal K}
\mapsto ((u_1, C_1 , z^1) , (u_2, C_2 , z^2 )) \in {\cal M}^{d_1}_{0 , k_{d_1}} (X)
\times {\cal M}^{d_2}_{0 , k_{d_2}} (X)$ and the morphism of evaluation
$ev^{d_1}_{k_{d_1}} \times ev^{d_2}_{k_{d_2}} : {\cal M}^{d_1}_{0 , k_{d_1}} (X)
\times {\cal M}^{d_2}_{0 , k_{d_2}} (X) \to X^{k_d}$, where $u_1 = u|_{C_1}$ and
$u_2 = u|_{C_2}$. It follows from the hypothesis that
the first morphism is injective. Also, from Lemma \ref{lemmacoker}, the morphism
$d|_{(u_2, C_2 , z^2 )}  ev^{d_2}_{k_{d_2}}$ is injective and $d|_{(u_1, C_1 , z^1)}
ev^{d_1}_{k_{d_1}}$ is surjective. We deduce from the hypothesis that the composition of these 
morphisms, that is $ev^d_{k_d}|_{\cal K}$, is injective. It suffices thus to prove that 
the restriction of $ev^d_{k_d}$ to a transversal of ${\cal K}$ at
the point $(u, C , z)$ is injective. Without loss of generality, we can assume that 
$z^1 = \{ z_1 , \dots , z_{k_{d_1}} \}$ and
$z^2 = \{ z_{k_{d_1} + 1} , \dots , z_{k_{d}} \}$.
Since $ev^d_{k_d} ({\cal K})$ is smooth in a neighborhood of $u(z)$ and its projection
on the first $k_{d_1}$ factors of $X^{k_d}$ is onto,
a transversal to this divisor can be chosen of the
form $\gamma : t \in \Delta \mapsto (x_1 (t), \dots , x_{k_d} (t)) \in X_{k_d}$ where only
one point $x_i (t)$, say $x_{k_d} (t)$, depends on $t$. The inverse image of this path is a 
smooth curve $B$ in 
$\overline{{\cal M}}^d_{0,k_d} (X)^*$ transversal to ${\cal K}$ at the point $(u, C , z)$.
Denote by $U \to B$ the restriction of $\overline{U}^d_{0 , k_d} (X) \to
\overline{{\cal M}}^d_{0,k_d} (X)^*$ over $B$ and by $s_1 , \dots , s_{k_d} : B \to U$
the tautological sections. The evaluation morphism $ev_U : U \to X$ contracts
the curves $Im (s_1) , \dots ,  Im (s_{k_d - 1})$ and we have to prove that it is injective
once restricted to $Im (s_{k_d})$. But the normal bundle of
$C_1$ (resp. $C_2$) in $U$ is isomorphic to ${\cal O}_{C_1} (-1)$ (resp. ${\cal O}_{C_2} (-1)$).
The morphism $dev_U$ maps this line bundle onto a subline bundle of degree at least
$k_{d_1} - 1$ (resp. $k_{d_2} - 2$) of $N_{u_1}$ (resp. $N_{u_2}$) since this morphism 
vanishes exactly at the points $z_1 , \dots , z_{k_{d_1}}$ (resp. $z_{k_{d_1} + 1} , \dots , 
z_{k_{d} - 1}$). Now this subline bundle contains the tangent line of $u (C_2)$ (resp. $u (C_1)$)
at the point $u (C_1) \cap u (C_2)$. Thus, it cannot be the unique subline bundle of degree
$k_{d_1}$ (resp. $k_{d_2} - 1$) of $N_{u_1}$ (resp. $N_{u_2}$). It follows that this subline 
bundle is of degree exactly $k_{d_1} -1$ (resp. $k_{d_2} - 2$) which implies that $d ev_U$ is
injective at every point of $C \setminus \{ z_1 , \dots , z_{k_d - 1} \}$. Hence the result. 
$\square$

\begin{rem}
With the help of an expression of the cokernel of the evaluation map at reducible curves
analogous to the one given in Lemma \ref{lemmacoker}, it would be possible to reduce the 
proofs of Propositions \ref{propcoker} and \ref{propcoker2} to some vanishing results, 
compare \S I.2 of \cite{Ko}.
\end{rem}

\subsection{Relation with the Gromov-Witten invariants}

\begin{prop}
\label{propGW}
Let $X$ be a convex manifold of dimension $3$ and $d \in H_2 (X ; \Z)$ be an effective 
homology class
such that $c_1 (X)d$ is even. Denote by $k_d = \frac{1}{2} c_1 (X)d$. Then, the following are
equivalents :

i) The morphism $ev^d_{k_d}$ is dominant.

ii) There exists an irreducible rational curve in the class $d$ which is balanced.

iii) The genus $0$ Gromov-Witten invariant $GW (X, d, pt, \dots , pt)$ does not vanish.
\end{prop}

{\bf Proof :}

{\it $i) \Rightarrow ii)$} Let $\underline{x} \in X^{k_d}$ be a regular value of $ev^d_{k_d}$ 
which does not
belong to the image of ${\cal K} = \overline{{\cal M}}^d_{0,k_d} (X) \setminus 
{\cal M}^d_{0,k_d} (X)^*$.
Then by hypothesis there exists $(u,C,z) \in {\cal M}^d_{0,k_d} (X)^*$ such that 
$u (z) = \underline{x}$
and $\coker d|_{(u,C,z)} ev^d_{k_d} = 0$. From lemma \ref{lemmacoker}, this means that
$H^1 (C , N_{u , -z}) = 0$. Now $N_{u , -z}$ is isomorphic to ${\cal O}_{C} (a - k_d) \oplus
{\cal O}_{C} (b - k_d)$ with $a + b \leq 2k_d - 2$. This condition thus implies that 
$a = b = k_d -1$,
which means that $u$ is immersed and balanced.

{\it $ii) \Rightarrow i)$} Follows from lemma \ref{lemmacoker} and the implicit function theorem.

{\it $iii) \Leftrightarrow i)$} Let $\underline{x} \in X^{k_d}$ be a regular value of $ev^d_{k_d}$ 
which does not belong to $ev^d_{k_d} ({\cal K})$. Then, by definition, $GW (X, d, pt, \dots , pt) =
\# (ev^d_{k_d})^{-1} (\underline{x})$. Hence the equivalence. $\square$\\

In particular, we deduce from Example $1$ of \S \ref{subsectbalanced} the following 
well known corollary. 

\begin{cor}
If $X = \C P^3$ and $d \neq 2[\C P^1] \in H_2 (X ; \Z)$ is effective, then the evaluation map
$ev^d_{k_d}$ is dominant and the genus $0$ Gromov-Witten invariant $GW (X, d, pt, \dots , pt)$ 
does not vanish. $\square$
\end{cor}

For some computations of such Gromov-Witten invariants, see \cite{RT}, \cite{FP} and 
references therein.

\section{Main results}

\subsection{Choice of a $Pin^-$ structure $\p$ on $\R X$}
\label{subsectspinstruct}

Let $(X , c_X)$ be a real algebraic convex $3$-manifold. Equip the real part $\R X$ with
a riemannian metric $g$ and denote by $R_X$ the $O_3 (\R)$-principal bundle of orthonormal
frames of $T \R X$.
Remember that the double covering of $O_3 (\R)$ which is non-trivial over each of its connected
components can be provided two different
group structures turning the covering map into a morphism. The one for which the lift of a 
reflexion is of order $4$ is denoted by $Pin^-_3$, see \cite{ABS}. The obstruction to lift the
$O_3 (\R)$-principal bundle $R_X$ into a $Pin^-_3$-principal bundle is given by the expression 
$w_2 (\R X) + w_1^2 (\R X) \in H^2 (\R X ; \Z /2\Z)$) where $w_1 (\R X)$ (resp. $w_2 (\R X)$)
stands for the first (resp. second) Stiefel-Whitney class of $\R X$,
see \cite{KirT} for instance. Now from Wu relations
(see \cite{Miln}, p. 132 for instance), the obstruction $w_2 (M) + w_1^2 (M)$ vanishes for every 
compact $3$-manifold $M$. From now on, we will thus fix such a $Pin^-_3$ structure $\p$ on 
$\R X$ and denote by $P_X$ the associated $Pin^-_3$-principal bundle. Note that on 
orientable components of $\R X$, the choice of an orientation allows to reduce the $Pin^-_3$ 
structure into a $Spin_3$ structure.

\subsection{Spinor states of balanced real rational curves}
\label{subsectspinorient} 

Let $(A , c_A) \subset (X , c_X)$ be a balanced irreducible real rational curve
realizing the homology class $d$ and with nonempty real part $\R A$. In particular,
$\R A$ is an immersed knot in $\R X$ and the normal bundle $N_A$ of $A$ is a
real holomorphic vector bundle isomorphic to ${\cal O}_A (k_d -1) \oplus 
{\cal O}_A (k_d -1)$ where $k_d = \frac{1}{2} c_1 (X)d \in \N^*$. 
Assume that $\R A$ has a marked point $x_0$ and fix an orientation on $T_{x_0} \R X$.
Since $w_1 (\R X)[\R A] = c_1 (X) [A] = 0 \mod (2)$, this orientation induces an orientation
on $T \R X|_{\R A}$.
The real ruled surface $P (N_A) \cong \C P^1 \times \C P^1$ has a real part homeomorphic
to a torus, it is thus isomorphic to $\C P^1 \times \C P^1$ equipped with 
the product of standard complex conjugation $conj \times conj$. Put an orientation on 
$\R A$, the real rank two vector bundle $\R N_A$ has an induced orientation such that
a direct orthonormal frame of $T_x \R A$ followed by a direct orthonormal frame of  
$\R N_A|_x$ provides a direct orthonormal frame of $T_x \R X$ for every $x \in \R A$.
Denote by $h$ the class of the real part of a section of $P (N_A)$ having vanishing 
self-intersection,
which is equipped with the orientation induced by the one of $\R A$. Similarly,
denote by $f$ the class of a real fiber of $P (\R N_A)$, which is equipped with the
orientation induced by the one of $\R N_A$. Then the couple $(h,f)$ provides a basis
of the lattice $H_1 (\R N_A ; \Z)$. Choose a real holomorphic subline bundle
$L \subset N_A$ such that $P(L) \subset P (N_A)$ is a section with vanishing self-intersection
if $k_d$ is odd and a section of bidegree $(1,1)$ whose real part is homologous
to $\pm (h+v) \in H_1 (\R N_A ; \Z)$ otherwise. The line bundle $L$ is then of degree
$k_d - 1$ if $k_d$ is odd and $k_d - 2$ otherwise. In both case, the real part $\R L$
is an orientable real line bundle. Denote by $\R E \supset T \R A$ the real rank two
sub-bundle of $T \R X$ which projects onto $\R L$. This bundle is equipped with a 
riemannian
metric induced by the one of $T \R X$. Choose then $(e_1 (p) , e_2 (p))_{p \in \R A}$
a loop of orthonormal frames of $\R E$ such that $(e_1 (p))_{p \in \R A}$ is a loop
of direct orthonormal frames of $T \R A$. Note that this choice is not unique since
$(e_1 (p) , - e_2 (p))_{p \in \R A}$ provides another such choice. Let 
$(e_3 (p))_{p \in \R A}$
be the unique section of $T \R X |_{\R A}$ such that 
$(e_1 (p) , e_2 (p) , e_3 (p))_{p \in \R A}$ is a loop
of direct orthonormal frames of $T \R X |_{\R A}$. Define then $sp (A) = +1$ if this
loop of $R_X$ lifts to a loop of the bundle $P_X$ and $sp (A) = -1$ otherwise.
This integer is called the {\it spinor state} of the curve $A$.
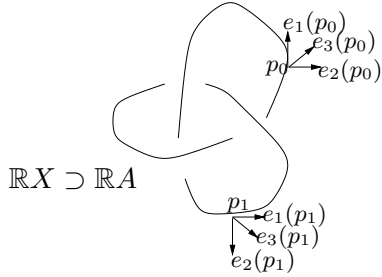
$$\vcenter{\hbox{\begin{picture}(0,0)%
\epsfig{file=conv1.pstex}%
\end{picture}%
\setlength{\unitlength}{1865sp}%
\begingroup\makeatletter\ifx\SetFigFont\undefined%
\gdef\SetFigFont#1#2#3#4#5{%
  \reset@font\fontsize{#1}{#2pt}%
  \fontfamily{#3}\fontseries{#4}\fontshape{#5}%
  \selectfont}%
\fi\endgroup%
\begin{picture}(4137,3605)(1996,-6172)
\put(5371,-5476){\makebox(0,0)[lb]{\smash{\SetFigFont{9}{10.8}{\rmdefault}{\mddefault}{\updefault}$e_1 (p_1)$}}}
\put(4906,-5296){\makebox(0,0)[lb]{\smash{\SetFigFont{9}{10.8}{\rmdefault}{\mddefault}{\updefault}$p_1$}}}
\put(1996,-5041){\makebox(0,0)[lb]{\smash{\SetFigFont{11}{13.2}{\rmdefault}{\mddefault}{\updefault}$\R X \supset \R A$}}}
\put(6031,-3181){\makebox(0,0)[lb]{\smash{\SetFigFont{9}{10.8}{\rmdefault}{\mddefault}{\updefault}$e_3 (p_0)$}}}
\put(5641,-2926){\makebox(0,0)[lb]{\smash{\SetFigFont{9}{10.8}{\rmdefault}{\mddefault}{\updefault}$e_1 (p_0)$}}}
\put(5416,-3481){\makebox(0,0)[lb]{\smash{\SetFigFont{9}{10.8}{\rmdefault}{\mddefault}{\updefault}$p_0$}}}
\put(4951,-6091){\makebox(0,0)[lb]{\smash{\SetFigFont{9}{10.8}{\rmdefault}{\mddefault}{\updefault}$e_2 (p_1)$}}}
\put(5281,-5761){\makebox(0,0)[lb]{\smash{\SetFigFont{9}{10.8}{\rmdefault}{\mddefault}{\updefault}$e_3 (p_1)$}}}
\put(6106,-3556){\makebox(0,0)[lb]{\smash{\SetFigFont{9}{10.8}{\rmdefault}{\mddefault}{\updefault}$e_2 (p_0)$}}}
\end{picture}
}}$$
It neither depends on the choice of $L$ nor on the choice of the orientation on
$\R A$ or on the metric $g$. It only depends on the pin structure $\p$ and 
on the choice of the orientation on $T_{x_0} \R X$ when $k_d$ is even. Indeed,
reversing this orientation change $sp (A)$ into its opposite in this case. Note that this
spinor state is inherited from the complex immersion of $A$ in $X$ and does depend on
this immersion, not just on the immersion of $\R A$ in $\R X$.

\begin{rem}
1) A spinor state for non-balanced real rational curve can be defined as well, but only
the balanced case will be relevant for our purpose.

2) A more geometric explanation of our construction can be given as follows. The choice
of a real holomorphic section $h \in P (N_A)$ with vanishing self-intersection provides
a canonical - homotopy class of - ribbon structure on our knot. In case our ribbon is a
M\oe bius strip, we have to cut this ribbon and perform half a twist in some direction
before gluing it again in order to get cylinder and be able to associate some framing.
This happens exactly when $k_d$ is even. The choice of the homology class
$\pm (h+v)$ instead of $\pm (h-v)$ in the above construction is equivalent to the choice
of such a direction. The point is that fixing an orientation on the tangent bundle over our
knot is enough to fix a choice of such a direction without ambiguity.
\end{rem}

\subsection{Statement of the results}

Let $(X , c_X)$ be a smooth real algebraic convex $3$-manifold. Denote by $(\R X)_1 , 
\dots , (\R X)_n$ the connected components of its real part $\R X$ and equip them with
a $Pin_3^-$ structure $\p$, see \S \ref{subsectspinstruct}. Let $d \in H_2 (X ; \Z)$ be
a class realized by real rational curves, such that $c_1 (X)d$ is even and $k_d = \frac{1}{2}
c_1 (X)d \in \N^*$.

Let $\underline{x}$ be a {\it real configuration} of $k_d$ distinct points of $X$, that is
a configuration of $k_d$ distinct points which are either real or exchanged by the involution 
$c_X$. For $i \in \{ 1 , \dots , n \}$, 
denote by $r_i = \# (\underline{x} \cap (\R X)_i)$ and assume that
$\sum_{i=1}^n r_i \neq 0$. Note that the $n$-tuple $r = (r_1 , \dots , r_n)$ encods
the equivariant isotopy class of $\underline{x}$. As soon as 
the choice of $\underline{x}$ is generic enough, there are only finitely many
connected rational curves in the homology class $d$ passing through 
$\underline{x}$. Moreover, these curves are all irreducible and balanced and their
number does not depend on the generic choice of $\underline{x}$, it is equal
to the genus $0$ Gromov-Witten invariant $N_d = GW (X,d, pt, \dots , pt)$.
Denote by ${\cal R}_d (\underline{x})$
the subset of these curves which are real and by  ${R}_d (\underline{x})$ its
cardinality. Since $\sum_{i=1}^n r_i \neq 0$ and $\underline{x}$ is generic, every curve 
$A \in {\cal R}_d (\underline{x})$ has nonempty real part. We can assume that
all the real points of
$\underline{x}$ are in a same connected component of $\R X$, say $(\R X)_1$, since otherwise
${\cal R}_d (\underline{x})$ is empty. 
Choose such a real 
point, say $x_{k_d}$, as well as an orientation of $T_{x_{k_d}} \R X$, and denote by 
$\underline{r}$ the equivariant isotopy class of $\underline{x}$ together with $x_{k_d}$ 
enriched with an orientation on $T_{x_{k_d}} \R X$.
Then, every curve $A \in {\cal R}_d (\underline{x})$ has a well defined
spinor state. Define finally :
$$\chi^{d,\p}_{\underline{r}} (\underline{x}) = \sum_{A \in {\cal R}_d (\underline{x})} sp(A) 
\in \Z.$$
\begin{theo}
\label{maintheorem}
Let $(X , c_X)$ be a smooth real algebraic convex $3$-manifold and $\p$ be a $Pin^-_3$ 
structure on $\R X$. Let $d \in H_2 (X ; \Z)$ be
such that $c_1 (X)d$ is even and different from $4$,
 and $k_d = \frac{1}{2} c_1 (X)d \in \N^*$. Let $\underline{x}= (x_1 , \dots , x_{k_d})$ be a real 
configuration of $k_d$ distinct points with at least one of them real, say $x_{k_d}$,
and $\underline{r}$ be the equivariant isotopy class of $\underline{x}$ together
with $x_{k_d}$ enriched with an orientation on $T_{x_{k_d}} \R X$
if $k_d$ is even. Then the integer 
$\chi^{d,\s}_{\underline{r}} (\underline{x})$ does not depend on the choice of 
$\underline{x}$, it only depends on $d, \p$ and $\underline{r}$.
The same holds when $c_1(X)d = 4$ and
the set of non-immersed rational curves of $X$ in the class $d$ is of codimension
at least $2$ in the moduli space ${\cal M}^d_{0,k_d} (X)$.
\end{theo}
(See Remark \ref{remarkdisc} for a discussion on the condition $c_1 (X)d \neq 4$. 
In fact,
I don't know any convex $3$-fold having a non-immersed curve whose homology class
$d$ satisfies $c_1 (X)d = 4$. Similarly, I don't know any real convex $3$-manifold having 
non-vanishing Gromov-Witten invariants and a non-connected real part.)\\

This integer $\chi^{d,\p}_r (\underline{x})$ is denoted by $\chi^{d,\p}_r$ and we set
$\chi^{d,\p}_r = 0$ as soon as it is not well defined.
\begin{cor}
\label{corvanish}
Under the assumptions of Theorem \ref{maintheorem}, assume in addition that $k_d$ is even and 
that $\R X$
has a non-orientable component $(\R X)_i$, $i \in \{1 , \dots , n\}$. If 
$r = (r_1 , \dots , r_n)$ with $r_i \neq 0$, then $\chi^{d,\p}_{\underline{r}} = 0$. In particular,
the genus $0$ Gromov-Witten invariant $GW(X,d,pt , \dots , pt)$ is even.
\end{cor}

{\bf Proof :}

Let $\underline{x} (t) = (x_1 (t), \dots , x_{k_d} (t))$, $t \in [0, 2 \pi]$, be a loop of 
configurations
of points of $X$ such that $x_2 (t), \dots , x_{k_d} (t)$ are fixed and $x_{k_d} (t)$ defines
a loop $\gamma$ of $\R X$ which is non-trivial against the first Stiefel-Whitney class of
$\R X$. Equip $T_{x_{k_d} (t)} \R X$ with an orientation depending continuously on $t$, then the
orientations of $T_{x_{k_d} (0)} \R X = T_{x_{k_d} (2 \pi)} \R X$ are opposite. Thus 
$\chi^{d,\p}_{\underline{r}} (\underline{x} (0)) = - \chi^{d,\p}_{\underline{r}} 
(\underline{x} (2 \pi))$. Now from
Theorem \ref{maintheorem}, $\chi^{d,\p}_{\underline{r}} (\underline{x} (0))  = 
 \chi^{d,\p}_{\underline{r}} (\underline{x} (2 \pi))$. This forces the vanishing 
$\chi^{d,\p}_{\underline{r}} =0$. $\square$\\

Denote by ${\cal R}_d^+ (\underline{x})$ (resp. ${\cal R}_d^- (\underline{x})$) the subset of
${\cal R}_d (\underline{x})$ consisting of curves having positive (resp. negative) spinor
state. Corollary \ref{corvanish} means that these two subsets have exactly same
cardinality, they are mirror to each other.
One example of such a manifold is the quadric in $\C P^4$ equipped with the real structure
whose real part is non-orientable. 

\begin{rem}
\label{centrem}
1) It is convenient to fix an orientation on orientable components of $\R X$ and to use it
to define the integer $\chi^{d,\p}_{\underline{r}}$. When $k_d$ is even, the sign of 
$\chi^{d,\p}_{\underline{r}}$ depends
on this choice of an orientation. We can then get rid of the necessity to put an orientation
on $T_{x_{k_d}} \R X$. Indeed, for orientable components, it is given by their orientation and for
non-orientable components, either $k_d$ is odd and the integer
$\chi^{d,\p}_{\underline{r}}$ is anyway well defined, or $k_d$ is even and
$\chi^{d,\p}_{\underline{r}} = 0$ from Corollary \ref{corvanish}.
Note that the pin structure $\p$ together with the orientation induce then a spin structure
on orientable components of $\R X$. {\bf From now on}, every orientable components of $\R X$ will
be equipped with a spin structure $\s$, which is more convenient for our need.
In case $\R X$ is orientable, the integer $\chi^{d,\p}_{\underline{r}}$ will be then rather 
denoted by $\chi^{d,\s}_r$.

2) Such a vanishing result as the one given by Corollary \ref{corvanish} can also be given 
for curves in $\C P^3$ of even degree. Indeed, the integers $\chi^{d,\s}_r (\underline{x}_0)$
and $\chi^{d,\s}_r (\underline{x}_1)$ must then have opposite signs if $\underline{x}_0$
and $\underline{x}_1$ are images to each other under a real reflexion of $\C P^3$. This has just
been noticed and communicated to me by G. Mikhalkin. Once more, this implies that the
associated Gromov-Witten invariant is even in even degree. The first degree where the values of
these invariants of the projective space are not known is thus $5$. The associated Gromov-Witten 
invariant is then $105$ (and $122129$ in degree $7$ as taken out from \cite{Gath}).
Note that similarly, the action of the group of real automorphisms of $X$ on the space of
real configuration of points provides symetries in the invariant $\chi^{\p}$.
\end{rem}

Thanks to this remark, we can define the
polynomial $\chi^{d,\p} (T) = \sum_{r \in \N^n} \chi^{d,\p}_r T^r \in \Z[T]$, where
$T^r = T_1^{r_1} \dots T_n^{r_n}$. This polynomial is of the same parity as the integer
$k_d$ and each of its monomials actually only depends on one indeterminate. Theorem
\ref{maintheorem} means that the function
$\chi^{\p} : d \in H_2 (X ; \Z) \mapsto \chi^{d,\p} (T) \in \Z[T]$ is an invariant
associated to the isomorphism class of the real algebraic convex $3$-manifold $(X , c_X)$.
As an application, this invariant provides the following lower bounds in real enumerative
 geometry.

\begin{cor}
\label{corlower}
Under the assumptions of Theorem \ref{maintheorem}, denote by ${R}_d (\underline{x})$
the number of real irreducible rational curves passing through $\underline{x}$ in the
class $d$ and by $N_d$ the associated Gromov-Witten invariant. Then,
$|\chi^{d,\p}_r| \leq R_d (\underline{x}) \leq N_d. \quad \square$
\end{cor}

Note that a similar invariant and similar lower bounds have already been obtained in 
\cite{Wels}, \cite{Wels2}  using real rational curves in real symplectic $4$-manifolds.
The question was then raised whether there exists such invariants in higher dimensions.
The results presented here thus provide a partial answer to this question.

The following natural questions arise from Corollary 
\ref{corlower}. Are the upper and lower bounds given by this corollary sharps ?
How to compute the invariant $\chi^{d,\p}_r$ ? See \cite{Sot} for a discussion of
related problems in real enumerative geometry and \cite{IKS} for an estimation
of the similar invariants constructed in \cite{Wels}, \cite{Wels2}. 

Finally, note that it is possible to 
understand the dependance of $\chi^{d,\p}_r$ with respect to $r$, see \S 
\ref{subsectfurth}.

\begin{rem}
G. Mikhalkin has just communicated to me that using considerations from tropical geometry,
he is able to prove that in degree $4$ in $\C P^3$, though the invariant $\chi^{4, \s}_8$
vanishes, the lower bound $0$ given by Corollary \ref{corlower} is sharp. This contrasts
with the complex dimension $2$, see \cite{IKS}.
\end{rem}

\section{Proof of Theorem \ref{maintheorem}}

Let $\tau \in {\cal S}_{k_d}$ having $\sum_{i=1}^n r_i$ fixed points in $\{ 1 , \dots , 
k_d \}$ and such that $\tau^2 = id$. Remember that $ev^d_{k_d} : 
(\overline{\cal M}^d_{0 , k_d} (X) , c_{\overline{\cal M} , \tau}) \to 
(X^{k_d} , c_\tau )$ is a real morphism between real algebraic varieties, see \S
\ref{subsectmoduli}. Denote by $\R_\tau X^{k_d}$ (resp. $\R_\tau 
\overline{\cal M}^d_{0 , k_d} (X) $) the real part of $(X^{k_d} , c_\tau )$ (resp.
$(\overline{\cal M}^d_{0 , k_d} (X) , c_{\overline{\cal M} , \tau})$) and by
$\R_\tau ev^d_{k_d} $ the restriction of $ev^d_{k_d}$ to 
$\R_\tau  \overline{\cal M}^d_{0 , k_d} (X) \to \R_\tau X^{k_d}$.

\subsection{Genericity arguments}
\label{subsectgeneric}

\begin{prop}
\label{propcusp}
Let $X$ be a smooth algebraic convex $3$-manifold and $u_0 : \C P^1 \to X$ be a morphism having
a unique cuspidal point at $z_0 \in \C P^1$. Assume that the holomorphic bundle $u_0^* TX \otimes
{\cal O}_{\C P^1} (-z_0)$ is generated by its global sections and denote by
$d = (u_0)_* [\C P^1]
\in H_2 (X ; \Z)$. Then, the locus of non-immersed curves is a subvariety of codimension $2$
of ${\cal M}^d_{0,0} (X)$ in the neighborhood of $u_0$.
\end{prop}

{\bf Proof :}

Denote by $T^*$ the holomorphic vector bundle of rank $3$ over ${\cal M}^d_{0,1} (X)$ whose
fiber over $(u , \C P^1 , z) \in {\cal M}^d_{0,1} (X)$ is the vector space $Hom (T_z \C P^1 , 
T_{u(z)} X)$. This bundle is equipped with a tautological section $\sigma : (u , \C P^1 , z)
\in {\cal M}^d_{0,1} (X) \mapsto d|_z u \in Hom_{\C} (T_z \C P^1 , 
T_{u(z)} X)$. The vanishing locus of $\sigma$ coincide with the locus of curves $(u , \C P^1 , z) 
\in {\cal M}^d_{0,1} (X)$ having a cuspidal point at $z$. Let us prove that $\sigma$ vanishes
transversely at the point $(u_0 , \C P^1 , z_0)$. For this purpose, we fix some holomorphic
local coordinates in the neighborhood of $z_0 \in \C P^1$ and in the neighborhood of $u_0 (z_0)
\in X$. These coordinates induce a connection $\nabla$ on the bundle $u_0^* TX$ in
the neighborhood of $z_0$, as well as a connection $\nabla^{T^*}$ on the bundle $T^*$
in the neighborhood of $(u_0 , \C P^1 , z_0)$. Let $\xi \in T_{z_0} \C P^1$ and $\zeta \in 
T_{u_0 (z_0)} X$. By hypothesis, there exists a section $v$ of the bundle $u_0^* TX$ such that
$v (z_0) = 0$ and $\nabla v|_{z_0} = \xi^* \otimes \zeta$. Then, $(v, \stackrel{.}{z} = 0) \in
T_{(u_0 , \C P^1 , z_0)} {\cal M}^d_{0,1} (X)$ and $\nabla^{T^*}_{(v, \stackrel{.}{z} = 0)}
\sigma = \nabla v|_{z_0} = \xi^* \otimes \zeta$. Hence $\nabla^{T^*} 
\sigma|_{(u_0 , \C P^1 , z_0)}$ is surjective, the result follows. $\square$

\begin{rem}
\label{remarkdisc}
1) Let $X=\C P^1 \times \C P^2$ and $d = kl$ where $k \geq 3$ and $l$ is the class of a line 
in a fiber $\{ z \} \times \C P^2$ of $X$. Then the locus of cuspidal curves is of 
codimension $1$
in ${\cal M}^d_{0,0} (X)$. However, if $(u , \C P^1)$ is such a curve and $z \in \C P^1$ is
a cuspidal point of $u$, then the bundle $u^* TX \otimes
{\cal O}_{\C P^1} (-z)$ is not generated by its global sections.\\

2)Let $(u , \C P^1 , z) \in {\cal M}^d_{0,k_d} (X)$ be such that $\dim H^1 (\C P^1 ; N_u
\otimes {\cal O}_{\C P^1} (-z)) = 1$. Then $u$ has a unique cuspidal point $z_0 \in \C P^1$
and $N_u \cong {\cal O}_{\C P^1} (k_d - 2) \oplus {\cal O}_{\C P^1} (k_d - 1)$. From the long
exact sequence associated to the short exact sequence
$0 \to T \C P^1  \to u^* TX \otimes {\cal O}_{\C P^1} (-z_0) \to 
N_u \otimes {\cal O}_{\C P^1} (-z_0) \to 0$ we deduce that the bundle $u^* TX \otimes 
{\cal O}_{\C P^1} (-z_0)$ is  generated by its global sections as soon as $k_d \geq 3$.
When $k_d = 1$, ${\cal M}^d_{0,k_d} (X)$ does not contain any non-immersed curve since
$\deg (u^* TX) = 2$, which implies that any morphism of sheaves $T \C P^1 \to u^* TX$ 
is injective. For $k_d = 2$, see the next remark.\\

3) Let $A$ be a cuspidal cubic curve in $\C P^2$. Denote by $Y$ the blown up of $\C P^2$
at $5$ distinct points of $A$ outside its cuspidal point and by $\widetilde{A}$ the strict 
transform of $A$ in $Y$. Then $\widetilde{A} \subset X = Y \times \C P^1$ satisfies
$c_1 (X) \widetilde{A} = 4$, but the
locus of non-immersed curves of $X$ in the class $[\widetilde{A}]$ is of codimension one.
In this case however, though $H^1 (\widetilde{A} ; TX|_{\widetilde{A}}) = 0$, $X$ is not convex.
\end{rem}

The divisor $\R_\tau {\cal K} = \R_\tau \overline{\cal M}^d_{0,k_d} (X) \setminus
\R_\tau {\cal M}^d_{0,k_d} (X)^*$ is real and made of reducible or multiple curves. Denote
by $\R_\tau {\cal K}_{reg}$ the locus of curves $(u , C , z) \in \R_\tau {\cal K}$ such that
$C$ has two irreducible components $C_1$ and $C_2$ for which $u(C_1)$, $u(C_2)$ are
real, immersed and meet in a single ordinary double point away from 
$u(z)$, and which have one of these two properties :

- Either $u(C_1)$, $u(C_2)$ are both 
balanced and for $i \in \{ 1,2\}$,
$z^i = C_i \cap z$ has cardinality $\frac{1}{2} c_1(X)d_i$ where $d_i = u_* [C_i] 
\in H_2 (X ; \Z)$. 

- Or $N_{u_1} \cong {\cal O}_{C_1} 
(k_{d_1} - 1) \oplus {\cal O}_{C_1} (k_{d_1})$ and $N_{u_2} \cong {\cal O}_{C_2} 
(k_{d_2} - 2) \oplus {\cal O}_{C_2} (k_{d_2} - 1)$ where $k_{d_1} = E( \frac{1}{2} c_1(X)d_1 )$,
$k_{d_2} = E( \frac{1}{2} c_1(X)d_2 ) + 1$, $u_1 = u|_{C_1}$ and $u_2 = u|_{C_2}$. Moreover,
$z^i = C_i \cap z$ has cardinality $k_{d_i}$ for $i \in \{ 1,2\}$, and the tangent line
to $u(C_1)$ (resp. to $u(C_2)$) at the point $u(C_1) \cap u(C_2)$ is not mapped to the unique
subline bundle of degree $k_{d_2} - 1$ (resp. $k_{d_1}$) of $N_{u_2}$ 
(resp. $N_{u_1}$). 

\begin{prop}
\label{propcodim2}
The image of the complement $\R_\tau {\cal K} \setminus \R_\tau {\cal K}_{reg}$ under $\R_\tau
ev^d_{k_d}$ is of codimension at least two in $\R_\tau X^{k_d}$.
\end{prop}

{\bf Proof :}

From Theorem \ref{theoKont}, ${\cal K}$ is a divisor with normal crossings of 
$\overline{\cal M}^d_{0,k_d} (X)$. As a consequence, the locus of curves having more than two
irreducible components or meeting in more than one single ordinary double point is of
codimension at least two in $\overline{\cal M}^d_{0,k_d} (X)$. It suffices then to prove the result
for curves $(u , C , z)$ such that $C$ has two irreducible components $C_1$ and $C_2$ for 
which $u(C_1)$, $u(C_2)$ meet in a single ordinary double point away from 
$u(z)$. For $i \in \{ 1,2\}$, denote by $z^i = C_i \cap z$, $d_i = u_* [C_i] 
\in H_2 (X ; \Z)$ and $k_{d_i} = \# z^i$. The morphism $ev^d_{k_d}$ restricted to these curves
is the composition of the restriction morphism $(u , C , z) \in {\cal K} \mapsto 
((u_1 , C_1 , z^1) , (u_2 , C_2 , z^2)) \in \overline{\cal M}^{d_1}_{0,k_{d_1}} (X) \times
\overline{\cal M}^{d_2}_{0,k_{d_2}} (X)$ and of the evaluation morphism $ev^{d_1}_{k_{d_1}} \times
ev^{d_2}_{k_{d_2}} : \overline{\cal M}^{d_1}_{0,k_{d_1}} (X) \times
\overline{\cal M}^{d_2}_{0,k_{d_2}} (X) \to X^{k_d}$. Two cases are now to be considered depending
on whether $c_1 (X)d_i$ is even or odd. In the first case, we can assume that $k_{d_i} =
\frac{1}{2} c_1 (X)d_i$. Indeed, if $k_{d_1} >
\frac{1}{2} c_1 (X)d_1$ for example, then $\dim \overline{\cal M}^{d_1}_{0,k_{d_1}} (X) =
c_1 (X)d_1 + k_{d_1} \leq 3k_{d_1} - 2 = \dim (X^{k_d})- 2$ so that already the image of
$ev^{d_1}_{k_{d_1}} \times
ev^{d_2}_{k_{d_2}}$ is of codimension $2$ in $X^{k_d}$. For the same reason, from Lemma 
\ref{lemmacoker}, we can assume that at least one of the two curves $u(C_1)$, $u(C_2)$, say
$u(C_1)$, is balanced. Note that perturbing the $k_{d_1}$ points $u (z^1)$ in $X^{k_{d_1}}$,
we can deform $u(C_1)$ in order to separate it from $u(C_2)$. It follows that the image of
the restriction morphism ${\cal K} \to \overline{\cal M}^{d_1}_{0,k_{d_1}} (X) \times
\overline{\cal M}^{d_2}_{0,k_{d_2}} (X)$ is of codimension at least one in 
$\overline{\cal M}^{d_1}_{0,k_{d_1}} (X) \times
\overline{\cal M}^{d_2}_{0,k_{d_2}} (X)$. Moreover, from Lemma \ref{lemmacoker}, the image of 
non-balanced curves of $\overline{\cal M}^{d_2}_{0,k_{d_2}} (X)$ under $ev^{d_2}_{k_{d_2}}$ is
of codimension at least $1$ in $X^{k_{d_2}}$. It follows from these two remarks that we can also
assume that $u(C_2)$ is balanced. The only remaining thing to prove in this first case is
that we can assume that both $u(C_1)$ and $u(C_2)$ are real. If this would not be the case, they 
would be exchanged by the involution $c_X$. They would thus meet  at each real point in the
configuration $\underline{x}$. By hypothesis, there exists such real points. Now since the
condition to have a marked point at the intersection point $C_1 \cap C_2$ is of codimension $2$
in $\overline{\cal M}^d_{0,k_{d}} (X)$ -it creates a ``ghost'' component-, we are done.

In the second case, we can for the same reason assume that $k_{d_1} = E( \frac{1}{2} c_1(X)d_1 )$,
$k_{d_2} = E( \frac{1}{2} c_1(X)d_2 ) + 1$ and then that $N_{u_1} \cong {\cal O}_{C_1} 
(k_{d_1} - 1) \oplus {\cal O}_{C_1} (k_{d_1})$ and $N_{u_2} \cong {\cal O}_{C_2} 
(k_{d_2} - 2) \oplus {\cal O}_{C_2} (k_{d_2} - 1)$. Moreover, we can assume that both
$u(C_1)$ and $u(C_2)$ are real. Finally, the condition on the tangent line
of $u(C_1)$ (resp. of $u(C_2)$) at the point $u(C_1) \cap u(C_2)$ can be obtained by
perturbation of the configuration of points $u(z^1)$ and $u(z^2)$ in $X^{k_{d_1}}$ and
$X^{k_{d_2}}$. It follows that the locus of curves which do not have all these properties
is of codimension at least two in $\overline{\cal M}^d_{0,k_{d}} (X)$. Hence the result. $\square$

\subsection{Proof of Theorem \ref{maintheorem}}

Let $\underline{x}^0 , \underline{x}^1 \in \R_\tau X^{k_d}$ be two regular values of
$\R_\tau ev^d_{k_d}$ which belong to the same connected component of $\R_\tau X^{k_d}$ 
outside the divisor $\R D_{reg} = \R_\tau ev^d_{k_d} (\R_\tau {\cal K}_{reg})$. 
We have to prove
that $\chi^{d , \p}_{\underline{r}} (\underline{x}^0) = \chi^{d , \p}_{\underline{r}} 
(\underline{x}^1)$. Let 
$\gamma : [0,1] \to \R_\tau X^{k_d}$ be a generic piecewise analytic path joining
$\underline{x}^0$ to $\underline{x}^1$ and transversal to $\R_\tau ev^d_{k_d}$. Denote by
$\R {\cal M}_\gamma$ the fiber product $\R_\tau \overline{\cal M}_{0,k_d}^d (X) \times_\gamma
[0,1]$ and by $\R \pi_\gamma$ the associated projection $\R {\cal M}_\gamma \to [0,1]$.
The integer $\chi^{d , \p}_{\underline{r}} (\gamma (t))$ is well defined for every 
$t \in [0,1]$ but a finite
number of parameters $0 < t_0 < \dots < t_j < 1$ corresponding either to critical values
of $\R \pi_\gamma$, or to the crossing of the wall $\R D_{reg} = \R_\tau ev^d_{k_d} 
(\R_\tau {\cal K}_{reg})$.
We can assume that the latter is crossed transversely by $\gamma$. Note also that
from Lemma \ref{lemmacoker} and Proposition \ref{propcusp}, critical points of $\R \pi_\gamma$
correspond to immersed irreducible curves $(u,C,z)$ with normal bundle isomorphic
to ${\cal O}_C (k_d - 2) \oplus {\cal O}_C (k_d)$.

\begin{lemma}
\label{lemmachoice}
The path $\gamma$ can be chosen such that when it crosses the wall of critical values
of $\R \pi_\gamma$, only one real point $x_i (t)$ of 
$\gamma (t) = (x_1 (t) , \dots , x_{k_d} (t))$ depends on $t$. Similarly, it can be chosen 
such that when it crosses the wall $\R D_{reg}$, either only one real point of 
$\gamma (t)$ depends on $t$, or only two points exchanged by $c_X$ depend on $t$. In the last case,
this choice can be done such that in addition to the $k_d - 2$ constant points of $\gamma$,
the corresponding family of curves in $\R {\cal M}_\gamma$ has a common fixed real point in 
$\R X$.
\end{lemma}

{\bf Proof :}

Let $(u , C , z)$, $C = C_1 \cup C_2$, be a stable map in $\R_\tau {\cal K}_{reg}$. The tangent 
space to $\R D_{reg}$ at the point $u(z)$ consists of infinitesimal deformations of $u(z)$
for which $u(C)$ deforms into a reducible connected curve passing through this configuration 
of points. From Proposition \ref{propcodim2}, two cases are then to be considered. If
$(u_1 , C_1 , z^1)$ and $(u_2 , C_2 , z^2)$ are balanced, where $z^i = z \cap C_i$
and $u_i = u|_{C_i}$,
we can assume that $u (C_2)$ has a real point in the configuration $u (z)$, say $u (z_{k_d})$.
There exists then a real section $v_2$ of the normal bundle $N_{u_2}$ which vanishes at every
point of $z^2$ except $z_{k_d}$ and which at the point $u(C_1) \cap u(C_2)$ does not belong to
the image of the tangent line of $u(C_1)$ in $N_{u_2}$. Then, the infinitesimal deformation of
$u(z_{k_d})$ in the direction $v_2 (z_{k_d})$, the other points $u(z_1) , \dots , u(z_{k_d - 1})$
being fixed, provides a vector transversal to $\R D_{reg}$ at the point $u(z)$. This proves
the lemma in this case.

Assume now that $N_{u_1} \cong {\cal O}_{C_1} (k_{d_1} - 1) \oplus {\cal O}_{C_1} (k_{d_1})$
and $N_{u_2} \cong {\cal O}_{C_2} (k_{d_2} - 2) \oplus {\cal O}_{C_2} (k_{d_2} - 1)$, where
$k_{d_i} = \# (z^i)$. If $u (z^2)$ has a real point, say $u (z_{k_d})$, then it suffices to
deform $u (z_{k_d})$ in a real direction which does not project onto a fiber of the unique
subline bundle of degree $k_{d_2} - 1$ of $N_{u_2}$, the other points 
$u(z_1) , \dots , u(z_{k_d - 1})$ being fixed, to define an infinitesimal deformation
of $u(z)$ transversal to $\R D_{reg}$. Otherwise, $u (z^2)$ has at least two points exchanged
by $c_X$, say $u (z_{{k_d} - 1})$ and $u (z_{k_d})$. Choose two distinct real points
$z'_{{k_d} - 1}$ and $z'_{k_d}$ in the real part of $C_2 \setminus z^2$, and denote by
$z' = (z_1 , \dots , z_{{k_d} - 2} , z'_{{k_d} - 1} , z'_{k_d})$. Then $(u , C , z') \in
\R_{\tau '} {\cal K}_{reg}$, where $\tau ' = \tau \circ ({k_d} - 1 \; {k_d})$. Let 
$(u_t , C_t , z'_t) \in \R_{\tau '} \overline{\cal M}^d_{0 , k_d} (X)$ be a path transversal to 
$\R_{\tau '} {\cal K}_{reg}$, such that only $u_t (z'_{k_d} (t))$ depends on $t$. Such a path exists
from what we have already done. Choose a deformation $z_{k_d - 1} (t)$, $z_{k_d} (t)$ of
the points $z_{k_d - 1}$, $z_{k_d}$ of $C_0 = C$ such that $u_t (z_{k_d - 1} (t))$
and $u_t (z_{k_d} (t))$ are exchanged by $c_X$. Then the path $(u_t , C_t , z (t)) \in 
\R_{\tau} \overline{\cal M}^d_{0 , k_d} (X)$, where $z(t) = (z_1 (t) , \dots , z_{k_d} (t))$,
is transversal to $\R_{\tau} {\cal K}_{reg}$ at the point $(u, C , z)$. Moreover, only the points
$u_t (z_{k_d - 1} (t))$ and $u_t (z_{k_d} (t))$ depend on $t$, and the curves $u_t (C_t)$ have an
additional real fixed point, namely $u_t (z'_{k_d - 1} (t))$. The lemma is thus proved
in this case.

Finally, assume that $(u , C , z)$ is a critical point of $\R \pi_\gamma$ and that
$u (z_{k_d})$ is real. Then $N_u \cong {\cal O}_{C} (k_d - 2) \oplus {\cal O}_{C} (k_d)$.
Choose any infinitesimal deformation of $u(z)$ which fixes $u(z_1) , \dots , u(z_{k_d - 1})$
and deforms $u (z_{k_d})$ in a real direction different from the one given by the unique
subline bundle of degree $k_d$ of $N_u$. This deformation is transversal to the wall
of critical values of $\R_\tau ev^d_{k_d}$. Hence the result. $\square$ \\

{\bf Proof of Theorem \ref{maintheorem} :}

Choose a path $\gamma$ given by Lemma \ref{lemmachoice}. We have to prove that the value
of the integer $\chi^{d , \p}_r (\gamma (t))$ does not change when $t$ crosses the parameters
$t_0 < \dots < t_j $. When this parameter corresponds to a critical value of $\R \pi_\gamma$,
this is given by the following Proposition \ref{propcrit}. When this parameter corresponds to
the crossing of the wall $\R_\tau D_{reg}$, this is given by the following Proposition 
\ref{propred}. $\square$

\begin{prop}
\label{propcrit}
Let $(u , C , z) \in \R {\cal M}_\gamma$ be a critical point of $\R \pi_\gamma$, $l_0 \in \N^*$
be the vanishing order of $d|_{(u , C , z)} \R \pi_\gamma$ and $t_0 = \R \pi_\gamma (u , C , z)
\in ] 0 , 1 [$. Then, there exists $\eta > 0$ and a neighborhood $W_0$ of $(u , C , z)$ in 
$\R {\cal M}_\gamma$  such that the following alternative occurs :

1) Either $l_0$ is odd, then for every $t \in ]t_0 - \eta , t_0 [$, $W_0 \cap \R \pi_\gamma^{-1}
(t) = \{(u_t^+ , C_t^+ , z_t^+) , (u_t^- , C_t^- , z_t^-) \}$ with $sp (u_t^+ (C_t^+)) =
-  sp (u_t^- (C_t^-))$ and for every $t \in ]t_0 ,  t_0 + \eta [$, $W_0 \cap \R \pi_\gamma^{-1}
(t) = \emptyset$, or the converse.

2) Or $l_0$ is even, then for every $t \in ]t_0 - \eta , t_0 + \eta [$, 
$W_0 \cap \R \pi_\gamma^{-1} (t) = \{(u_t , C_t , z_t) \}$ and $sp (u_t (C_t))$ does not
depend on $t \neq t_0$.
\end{prop}

{\bf Proof :}

Let us complexify the path $\gamma$ in the neighborhood of $t_0 \in ]0,1[$ in order to get
an analytic path $\gamma_\C : t \in \Delta_{t_0} (\eta) = \{ z \in \C \, | \, |z - t_0| < \eta \}
\to X^{k_d}$ such that for every $t \in \Delta_{t_0} (\eta)$, $\gamma_\C (\overline{t}) =
c_\tau \circ \gamma_\C (t)$ and $\gamma_\C |_{]t_0 - \eta , t_0 + \eta [} = \gamma$. Denote
then by ${\cal M}_\gamma$ the fiber product $\overline{\cal M}_{0,k_d}^d (X) \times_\gamma
[0,1]$ and by $\pi_\gamma$ the associated projection ${\cal M}_\gamma \to \Delta_{t_0} (\eta)$.
Denote also by $U \to {\cal M}_\gamma$ the restriction of $\overline{U}_{0,k_d}^d (X) \to
\overline{\cal M}_{0,k_d}^d (X)^* $. Hence, $U$ is a ruled surface over the smooth curve 
${\cal M}_\gamma$. Denote by $c_\gamma$ and $c_U$ the real structures on ${\cal M}_\gamma$
and $U$, so that the submersion $U \to {\cal M}_\gamma$ is $\Z / 2\Z$-equivariant. Fix a real
parametrization $\lambda \in \Delta \mapsto (u_\lambda , C_\lambda , z^\lambda) \in
{\cal M}_\gamma$ so that the projection ${\cal M}_\gamma \to \Delta_{t_0} (\eta)$ writes
$\lambda \in \Delta \mapsto \eta \lambda^{l_0 + 1} + t_0$. Denote then by $N$ the holomorphic
rank two vector bundle over $U$ whose fiber over $C_\lambda$ is the normal bundle $N_{u_\lambda}$.
For $\lambda \neq 0$, we have $N_{u_\lambda} \cong {\cal O}_{C_\lambda} (k_d - 1) \oplus
{\cal O}_{C_\lambda} (k_d - 1)$ and for $\lambda = 0$, $N_{u_0} \cong 
{\cal O}_{C_0} (k_d - 2) \oplus {\cal O}_{C_0} (k_d)$. Thus, the projectivization $P(N)$
is a deformation of ruled surfaces over the disk $\Delta$ such that the fibers
$P(N_{u_\lambda})$ over $\lambda \in \Delta \setminus \{ 0 \}$ are isomorphic to $\C P^1
\times \C P^1$ and the fiber $P(N_{u_0})$ is isomorphic to  the ruled surface
$P({\cal O}_{\C P^1} \oplus {\cal O}_{\C P^1} (2))$. The path $\gamma_\C$ can be written
$(x_1 (t) , \dots , x_{k_d} (t))$ where only the last point $x_{k_d} (t) \in \R X$ depends on
$t$ and $x'_{k_d} (t)$ is never tangent to a rational curve of ${\cal M}_\gamma$. It follows
that each ruled surface $P(N_{u_\lambda})$ has a marked point $w_\lambda$ which is the
projectivization of the complex line generated by $x'_{k_d} (t_\lambda) \in N_{u_\lambda}$,
for $t_\lambda = \eta \lambda^{l_0 + 1} + t_0$. Since $\gamma_\C$ is transversal to the
wall of critical values of $ev^d_{k_d}$, $w_0$ does not belong to the exceptional section
of $P(N_{u_0})$. Let $h_0$ be a real section of $P(N_{u_0})$ disjoint from the exceptional 
section and passing through $w_0$. This section can be deformed into an analytic family
$(h_\lambda)_{\lambda \in \Delta}$ of sections of $P(N_{u_\lambda})$ such that $h_\lambda$
passes through $w_\lambda$, $h_{\overline{\lambda}} = \overline{h_\lambda}$ and $h_\lambda$
is of bidegree $(1,1)$ in $P(N_{u_\lambda})$ for $\lambda \neq 0$. Denote by $L_\lambda \subset
N_{u_\lambda}$ the subline bundle associated to $h_\lambda$ and by $M_\lambda =
N_{u_\lambda} / L_\lambda$. We have $\deg (L_\lambda) = k_d -2$ and $\deg (M_\lambda) = k_d$.
Moreover, the extension $0 \to L_\lambda \to N_{u_\lambda} \to M_\lambda \to 0$ is
non-trivial for $\lambda \neq 0$ and trivial for $\lambda = 0$. Fix an identification
$U \cong \Delta \times \C P^1$ such that $z^\lambda$ is constant and let $\delta$ be a real generator
of $H^1 (\C P^1 ; {\cal O}_{\C P^1} (-2))$. There exists a holomorphic function $f : \Delta \to
\C$ such that the extension class of $0 \to L_\lambda \to N_{u_\lambda} \to M_\lambda \to 0$ is
$f (\lambda) \delta \in H^1 (\C P^1 ; M_\lambda^* \otimes L_\lambda) = 
H^1 (\C P^1 ; {\cal O}_{\C P^1} (-2))$. Then, $f(0) = 0$, $f(\lambda) \neq 0$ if $\lambda \neq 0$
and we define $l_f \in \N^*$ to be the vanishing order of $f$ at $0$. Let $V_0 , V_1$ be the
two standard affine charts of $\C P^1$, chosen such that $z_{k_d} (t_0)$ belongs
to $V_0$. The deformation $N$ is then obtained as the gluing of the trivializations
${\cal V}_0 = \Delta \times V_0 \times {\cal O}_{\C P^1} (k_d - 2) \oplus {\cal O}_{\C P^1} (k_d)$
and ${\cal V}_1 = \Delta \times V_1 \times {\cal O}_{\C P^1} (k_d - 2) \oplus 
{\cal O}_{\C P^1} (k_d)$ with gluing maps :
$$(\lambda , z , (s_1 , s_2)) \in {\cal V}_0 \cap {\cal V}_1 \mapsto \big( \lambda , z ,
(s_1 , s_1 + f(t) \delta s_2)\big) \in {\cal V}_1 \cap {\cal V}_0.$$
Now, note that the bundle $N$ has a tautological section 
$v^\lambda = \frac{d}{d\lambda} u_\lambda$. This section vanishes at the fixed points $z^\lambda_1 , 
\dots , z^\lambda_{k_d - 1}$. Denote its coordinates in the trivialization ${\cal V}_0$
by $v^\lambda (z) = (\lambda , z , (s^\lambda_1 (z) , s^\lambda_2 (z)))$. Since $(u_0 , C_0 ,
z^0)$ is a critical point of $\pi_\gamma$, from Lemma \ref{lemmacoker}, $v_0$ defines a non-zero
element of $H^0 (C_0 ; N_{u_0} \otimes {\cal O}_{C_0} (-z_0)) = H^0 (\C P^1 ; {\cal O}_{\C P^1}
(-2) \oplus {\cal O}_{\C P^1})$. It follows that the component $s^0_1$ of $v^0$ vanishes. Let us
prove that the vanishing order of $s^\lambda_1 (z_{k_d})$ at $\lambda = 0$ is $l_f$. Indeed, in
the Taylor expansion of $s^\lambda_1$ all the terms of order less than $l_f$ define sections
of $L_\lambda$ which vanish at the points $z^\lambda_1 , 
\dots , z^\lambda_{k_d - 1}$. These are thus the zero section. Denote by $(s_1)^{(l_f)}$
(resp. $(s_2)^{(0)}$) the term of order $l_f$ (resp. $0$) in the Taylor expansion of $s^\lambda_1$
(resp. $s^\lambda_2$). Then the triple $(\lambda , z , ((s_1)^{(l_f)} (z^\lambda_{k_d}) ,
(s_2)^{(0)} (z^\lambda_{k_d}))$ defines a section of the bundle $N$. If $(s_1)^{(l_f)} (z_{k_d}) = 0$,
then this section vanishes at the $k_d$ points $z_1 , \dots , z_{k_d}$. Since when $\lambda
\neq 0$, $N_{u_\lambda} \cong {\cal O}_{C_\lambda} (k_d - 1) \oplus {\cal O}_{C_\lambda} (k_d - 1)$,
this is impossible and thus the vanishing order of $s^\lambda_1 (z_{k_d})$ at $\lambda = 0$ is 
exactly $l_f$. Note that since by construction $L_\lambda$ contains $x'_{k_d} (t_\lambda)$,
where $t_\lambda = \eta \lambda^{l_0 + 1} + t_0$, the second component $s^\lambda_2$ is
identically zero at the point $z^\lambda_{k_d}$. From $\pi_\gamma (u_\lambda , C_\lambda , z^\lambda)
= u_\lambda (z^\lambda)$, we deduce that $d|_{(u_\lambda , C_\lambda , z^\lambda)} \pi_\gamma
(v_\lambda) = v_\lambda (z^\lambda) = (0 , \dots, 0 , v_\lambda (z^\lambda_{k_d}))$. Hence the
vanishing order $l_0$ of $d|_{(u , C , z)} \pi_\gamma$ is the one of $s^\lambda_1 (z^\lambda_{k_d})$
at the point $\lambda = 0$, that is $l_f$. 

Now, fix an orientation of the curves $\R C_\lambda$, $\lambda \in ]-1,1[$, and denote by 
$[\R v_\lambda]$ the section of $P (\R N_{u_\lambda})$ associated to the line bundle generated
by $v_\lambda$. The latter is equipped with the orientation induced by the one of $\R C_\lambda$.
Note that it is the real part of a section  of $P (N_{u_\lambda})$ with self-intersection zero,
since the holomorphic line bundle generated by $v_\lambda$ has degree $k_d - 1$. 
Similarly, the real rank
two vector bundle $\R N_{u_\lambda}$ has an  orientation induced by the ones of $\R C_\lambda$
and $T_{x_{k_d}} \R X$. Let $f_\lambda$ be a fiber of $P (\R N_{u_\lambda})$
equipped with the orientation induced by the one of $\R N_{u_\lambda}$. Then, for $\lambda \neq 0$,
$([\R v_\lambda] , f_\lambda)$ provides a basis of the lattice $H_1 (P (\R N_{u_\lambda}) ; \Z)$.
The section $\R h_\lambda = P (\R L_\lambda)$ has bidegree $(1,1)$. Let $\epsilon \in \{ \pm 1 \}$
be such that for $\lambda < 0$, $[\R h_\lambda] = [\R v_\lambda] + \epsilon  f_\lambda \in
H_1 (P (\R N_{u_\lambda}) ; \Z)$. Then, for $\lambda > 0$, 
$[\R h_\lambda] = [\R v_\lambda] + (-1)^{l_f} \epsilon  f_\lambda \in 
H_1 (P (\R N_{u_\lambda}) ; \Z)$. This follows from the local expression of $v_\lambda$ in the
neighborhood of $z^\lambda_{k_d}$, which is of the form $z \in \R C_\lambda \mapsto
(\lambda^{l_f} , z)$, whereas the classes $[\R h_\lambda]$ and $f_\lambda$ continuously depend
on $\lambda \in ]-1 , 1[$.
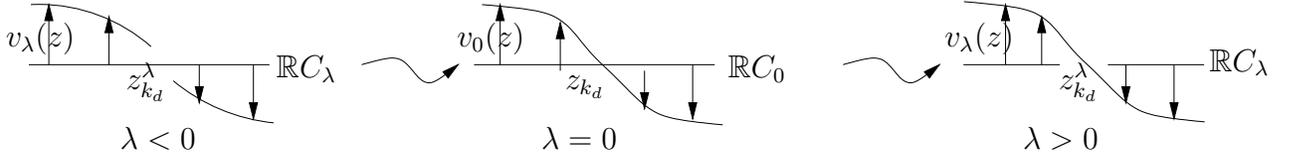
$$\vcenter{\hbox{\begin{picture}(0,0)%
\epsfig{file=conv2.pstex}%
\end{picture}%
\setlength{\unitlength}{3315sp}%
\begingroup\makeatletter\ifx\SetFigFont\undefined%
\gdef\SetFigFont#1#2#3#4#5{%
  \reset@font\fontsize{#1}{#2pt}%
  \fontfamily{#3}\fontseries{#4}\fontshape{#5}%
  \selectfont}%
\fi\endgroup%
\begin{picture}(9000,1114)(1396,-5191)
\put(1396,-4426){\makebox(0,0)[lb]{\smash{\SetFigFont{12}{14.4}{\rmdefault}{\mddefault}{\updefault}$v_\lambda (z)$}}}
\put(4771,-4426){\makebox(0,0)[lb]{\smash{\SetFigFont{12}{14.4}{\rmdefault}{\mddefault}{\updefault}$v_0 (z)$}}}
\put(6796,-4651){\makebox(0,0)[lb]{\smash{\SetFigFont{12}{14.4}{\rmdefault}{\mddefault}{\updefault}$\R C_0$}}}
\put(3421,-4651){\makebox(0,0)[lb]{\smash{\SetFigFont{12}{14.4}{\rmdefault}{\mddefault}{\updefault}$\R C_\lambda$}}}
\put(5581,-4741){\makebox(0,0)[lb]{\smash{\SetFigFont{12}{14.4}{\rmdefault}{\mddefault}{\updefault}$z_{k_d}$}}}
\put(10396,-4606){\makebox(0,0)[lb]{\smash{\SetFigFont{12}{14.4}{\rmdefault}{\mddefault}{\updefault}$\R C_\lambda$}}}
\put(9271,-4741){\makebox(0,0)[lb]{\smash{\SetFigFont{12}{14.4}{\rmdefault}{\mddefault}{\updefault}$z^{\lambda}_{k_d}$}}}
\put(2296,-4741){\makebox(0,0)[lb]{\smash{\SetFigFont{12}{14.4}{\rmdefault}{\mddefault}{\updefault}$z^{\lambda}_{k_d}$}}}
\put(8416,-4426){\makebox(0,0)[lb]{\smash{\SetFigFont{12}{14.4}{\rmdefault}{\mddefault}{\updefault}$v_\lambda (z)$}}}
\put(9001,-5191){\makebox(0,0)[lb]{\smash{\SetFigFont{12}{14.4}{\rmdefault}{\mddefault}{\updefault}$\lambda > 0$}}}
\put(5401,-5191){\makebox(0,0)[lb]{\smash{\SetFigFont{12}{14.4}{\rmdefault}{\mddefault}{\updefault}$\lambda = 0$}}}
\put(2251,-5191){\makebox(0,0)[lb]{\smash{\SetFigFont{12}{14.4}{\rmdefault}{\mddefault}{\updefault}$\lambda < 0$}}}
\end{picture}
}}$$
It follows that if $l_0 = l_f$ is even, the spinor states of 
$(u_\lambda , C_\lambda , z^\lambda)$ and $(u_{- \lambda} , C_{- \lambda} , z^{- \lambda})$
coincide for every $\lambda \neq 0$, whereas they are opposite when $l_0 = l_f$ is odd. Indeed,
in this last case, the homotopy classes of the loops of the bundle $R_X$ of orthonormal 
frames that we construct
in order to compute the spinor state (see \S \ref{subsectspinorient}) exactly differ from
the class of a non-trivial loop in a fiber of $R_X$. The result now follows from the fact that
the projection $\pi_\gamma$ has been identified in the neighborhood of 
$(u , C , z) \in \R {\cal M}_\gamma$ with the map 
$\lambda \in ]-1 , 1 [ \mapsto \eta \lambda^{l_0 + 1} + t_0 \in ]t_0 - \eta , t_0 + \eta [$.
$\square$

\begin{rem}
From the proof of Proposition \ref{propcrit} arise the following question. Are the generic
critical points of the evaluation map $ev^d_{k_d} : \overline{\cal M}^d_{0, k_d} (X) \to X^{k_d}$
non-degenerate when $X$ is convex, e.g. when $X = \C P^3$ ?
\end{rem}

\begin{prop}
\label{propred}
Let $(u , C , z) \in \R {\cal M}_\gamma$ be such that $C$ is reducible and 
$t_1 = \R \pi_\gamma (u , C , z) \in ] 0 , 1 [$. Then, there exists $\eta > 0$ and a 
neighborhood $W_1$ of $(u , C , z)$ in $\R {\cal M}_\gamma$  such that
for every $t \in ]t_1 - \eta , t_1 + \eta [$, $W_1 \cap \R \pi_\gamma^{-1}
(t) = \{(u_t , C_t , z_t)  \}$. Moreover, $sp (u_t (C_t))$ does not depend on $t \neq t_1$.
\end{prop}

{\bf Proof :}

From Propositions \ref{propcoker}, \ref{propcoker2} and the choice of $\gamma$, $(u , C , z)$ is
a regular point of $\R \pi_\gamma$. Thus there exists $\eta > 0$ and a 
neighborhood $W_1$ of $(u , C , z)$ in $\R {\cal M}_\gamma$  such that
for every $t \in ]t_1 - \eta , t_1 + \eta [$, $W_1 \cap \R \pi_\gamma^{-1}
(t) = \{(u_t , C_t , z_t)  \}$. Denote by $U \to W_1$ the universal curve, it has a unique 
singular 
fiber over $t = t_1$. The real part $\R U$ is thus homeomorphic to a cylinder $Cyl$ blown up
at one point, that is to the connected sum of a cylinder and a Moebius strip. From Lemma
\ref{lemmachoice}, the family $u_t (C_t)$ passes through $k_d - 1$ fixed points. For every 
$t \neq t_1$, the restriction $T \R U |_{\R C_t}$ projects onto a subline bundle of the normal 
bundle
of $\R C_t$ in $\R X$ which is the real part of a holomorphic line bundle of degree $k_d - 1$.
Assume that $k_d$ is odd, so that this line bundle can be used in order to compute the 
spinor state of $\R C_t$, $t \neq t_1$. The case $k_d$ even will follow along the same lines.
Fix a continuous family $(e_1^t (z) , e_2^t (z))$,
$t \neq t_1$, $z \in \R C_t$, of orthonormal frames of the tangent space at $u_t (z)$ of the image
of
$\R U$ in $\R X$. This family is completed in a family $(e_1^t (z) , e_2^t (z), e_3^t (z))$
of orthonormal frames of $T_{u_t (z)} \R X$. Remove a 
small disk $D$
around the blown up point of $Cyl$. In $H_1 (Cyl \setminus D ; \Z / 2\Z)$, the class of the
fiber $[\R C_{t_1 + \eta}]$ equals the sum of the class of the
fiber $[\R C_{t_1 - \eta}]$ and the class of the boundary $\partial D$. Equip these loops with
a family of orthonormal frames of $T \R X$ as before and denote
by $p (\R C_t)$ and $p (\partial D)$ the loops of $R_X$ thus obtained. These loops are related
by 
$p (\R C_{t_1 + \eta}) = p (\R C_{t_1 - \eta}) + p (\partial D) + p_0 
\in H_1 (R_X ; \Z / 2\Z),$ where $p_0$ is a non-trivial loop in a fiber of $R_X$. This is suggested
by the following picture.
$$\vcenter
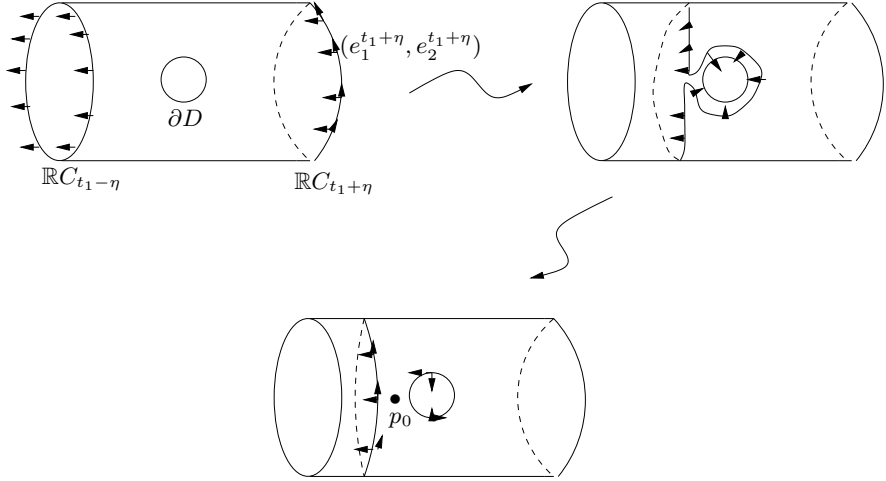
{\hbox{\begin{picture}(0,0)%
\epsfig{file=conv3.pstex}%
\end{picture}%
\setlength{\unitlength}{2486sp}%
\begingroup\makeatletter\ifx\SetFigFont\undefined%
\gdef\SetFigFont#1#2#3#4#5{%
  \reset@font\fontsize{#1}{#2pt}%
  \fontfamily{#3}\fontseries{#4}\fontshape{#5}%
  \selectfont}%
\fi\endgroup%
\begin{picture}(8751,4749)(2059,-8173)
\put(5896,-7621){\makebox(0,0)[lb]{\smash{\SetFigFont{9}{10.8}{\rmdefault}{\mddefault}{\updefault}$p_0$}}}
\put(5896,-7441){\makebox(0,0)[lb]{\smash{\SetFigFont{9}{10.8}{\rmdefault}{\mddefault}{\updefault}$\bullet$}}}
\put(2431,-5236){\makebox(0,0)[lb]{\smash{\SetFigFont{9}{10.8}{\rmdefault}{\mddefault}{\updefault}$\R C_{t_1 - \eta}$}}}
\put(4951,-5281){\makebox(0,0)[lb]{\smash{\SetFigFont{9}{10.8}{\rmdefault}{\mddefault}{\updefault}$\R C_{t_1 +\eta}$}}}
\put(3691,-4651){\makebox(0,0)[lb]{\smash{\SetFigFont{9}{10.8}{\rmdefault}{\mddefault}{\updefault}$\partial D$}}}
\put(5401,-3931){\makebox(0,0)[lb]{\smash{\SetFigFont{9}{10.8}{\rmdefault}{\mddefault}{\updefault}$(e_1^{t_1+\eta} , e_2^{t_1 +\eta})$}}}
\end{picture}
}}$$
Assume now that $\R X$ is orientable, or that if $C_1 , C_2$ denote the 
two components of $C$, then the curves $(u|_{C_1} , C_1 , z \cap C_1)$ and 
$(u|_{C_2} , C_2 , z \cap C_2)$ are balanced. In the Moebius strip that is glued instead
of the disk $D$, the boundary $\partial D$ is homologous to two times the core $A$ of the Moebius
strip. Denote by $p (A)$ a loop in $R_X$ obtained by equipping every point $z \in A$ with 
a orthonormal frame $(e_1 (z) , e_2 (z), e_3 (z))$ of $T_{u_t (z)} \R X$,
such that $e_1 (z) \in T_z A$ and $(e_1 (z) , e_2 (z))$ generates the tangent space of the Moebius
strip at $z$ except in a neighborhood of a point $z_0 \in A$ where $(e_1 (z) , e_2 (z), e_3 (z))$
does a half twist around the axis $TA$ (which is necessary since the Moebius strip is not
orientable). Then, one observes the relation
$p (\partial D) = 2p(A) + p_0 \in H_1 (R_X ; \Z / 2\Z).$
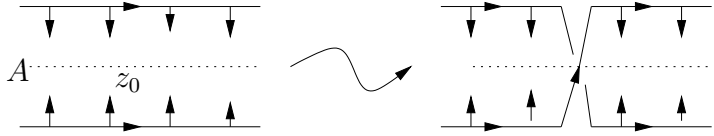
$$\vcenter{\hbox{\begin{picture}(0,0)%
\epsfig{file=conv4.pstex}%
\end{picture}%
\setlength{\unitlength}{3315sp}%
\begingroup\makeatletter\ifx\SetFigFont\undefined%
\gdef\SetFigFont#1#2#3#4#5{%
  \reset@font\fontsize{#1}{#2pt}%
  \fontfamily{#3}\fontseries{#4}\fontshape{#5}%
  \selectfont}%
\fi\endgroup%
\begin{picture}(5277,924)(1936,-5473)
\put(1936,-5101){\makebox(0,0)[lb]{\smash{\SetFigFont{12}{14.4}{\rmdefault}{\mddefault}{\updefault}$A$}}}
\put(2746,-5146){\makebox(0,0)[lb]{\smash{\SetFigFont{12}{14.4}{\rmdefault}{\mddefault}{\updefault}$z_0$}}}
\end{picture}
}}$$
We finally deduce that $p (\R C_{t_1 + \eta}) = p (\R C_{t_1 - \eta}) \in H_1 (R_X ; \Z / 2\Z),$
so that $sp(u_{t_1 + \eta} (C_{t_1 + \eta})) = sp(u_{t_1 - \eta} (C_{t_1 - \eta}))$.

It remains to consider the case when $\R X$ is not orientable and the curves 
$(u|_{C_1} , C_1 , z \cap C_1)$ and 
$(u|_{C_2} , C_2 , z \cap C_2)$ are not balanced. In this case, the double cover of
orientations $\widetilde{\R X}
\to \R X$ is non-trivial over $A_1 = u(\R C_1)$ and $A_2 = u(\R C_2)$. Denote by 
$\widetilde{A}_1$ and $\widetilde{A}_2$ the lifts of $A_1$ and $A_2$ in $\widetilde{\R X}$.
These are two immersed circles in $\widetilde{\R X}$ which intersects in two ordinary double points
$a_1$ and $a_2$. These points divide $\widetilde{A}_1$ (resp. $\widetilde{A}_2$) in two
connected components $\widetilde{A}_1^+$ and $\widetilde{A}_1^-$ (resp. $\widetilde{A}_2^+$ and 
$\widetilde{A}_2^-$). Let $\tilde{x}_{k_d}$ be the lift of $x_{k_d}$ in $\widetilde{\R X}$ given by
the choice of an orientation on $T_{x_{k_d}} \R X$. For $t \neq t_1$ in the neighborhood of $t_1$, 
denote by $\widetilde{\R C}_t$
the unique lift of $u (\R C_t)$ in $\widetilde{\R X}$ which passes through $\tilde{x}_{k_d}$. 
Without loss of generality, we can assume that $\tilde{x}_{k_d} \in
\widetilde{A}_1^+$ and that for $t_+ \in ]t_1 , t_1 + \eta [$ (resp. $t_- \in ]t_1 - \eta, t_1 [$),
$\widetilde{\R C}_{t_+}$ is in the neighborhood of $\widetilde{A}_1^+ \cup \widetilde{A}_2^+$
(resp. $\widetilde{A}_1^+ \cup \widetilde{A}_2^-$). Denote by $p (\widetilde{\R C}_{t_\pm})$
the loop of $R_X$ obtained by equipping the curve $\widetilde{\R C}_{t_\pm}$ with a family of
direct orthonormal frames $(e_1 (z) , e_2 (z), e_3 (z))_{z \in \widetilde{\R C}_{t_\pm}}$ as before.
Then, the difference $\lim_{t_+ \to t_1} p (\widetilde{\R C}_{t_+}) -
\lim_{t_- \to t_1} p (\widetilde{\R C}_{t_-})$ is observed to be a loop of $R_X$ over
$\widetilde{A}_2$ which is invariant under the involution $inv_R : (e_1 , e_2 , e_3) \in R_X|_x 
\mapsto (d|_x \inv (e_1) , d|_x \inv (e_2) , - d|_x \inv (e_3)) \in R_X|_{inv(x)} $ where
$inv$ is the involution of the covering $\widetilde{\R X} \to \R X$. Indeed, its restrictions
$p (\widetilde{A}_2^+)$ and $p (\widetilde{A}_2^-)$ over $\widetilde{A}_2^+$ and 
$\widetilde{A}_2^-$ respectively are paths exchanged by $inv_R$.
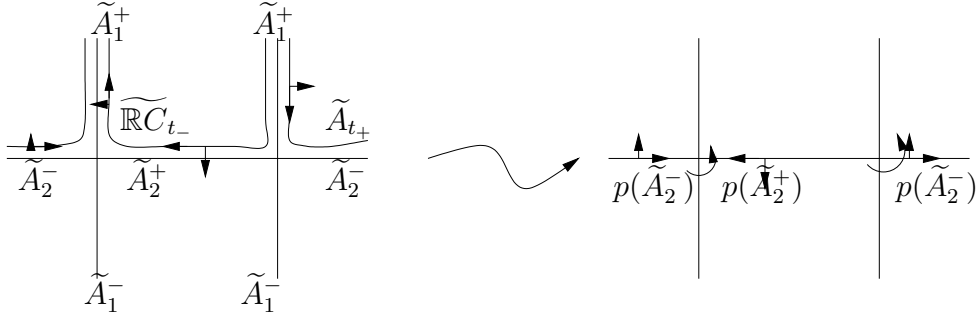
$$\vcenter{\hbox{\begin{picture}(0,0)%
\epsfig{file=conv5.pstex}%
\end{picture}%
\setlength{\unitlength}{3315sp}%
\begingroup\makeatletter\ifx\SetFigFont\undefined%
\gdef\SetFigFont#1#2#3#4#5{%
  \reset@font\fontsize{#1}{#2pt}%
  \fontfamily{#3}\fontseries{#4}\fontshape{#5}%
  \selectfont}%
\fi\endgroup%
\begin{picture}(7224,2250)(2014,-5641)
\put(2611,-5641){\makebox(0,0)[lb]{\smash{\SetFigFont{12}{14.4}{\rmdefault}{\mddefault}{\updefault}$\widetilde{A}^-_1$}}}
\put(3781,-5641){\makebox(0,0)[lb]{\smash{\SetFigFont{12}{14.4}{\rmdefault}{\mddefault}{\updefault}$\widetilde{A}^-_1$}}}
\put(4411,-4786){\makebox(0,0)[lb]{\smash{\SetFigFont{12}{14.4}{\rmdefault}{\mddefault}{\updefault}$\widetilde{A}^-_2$}}}
\put(2656,-3616){\makebox(0,0)[lb]{\smash{\SetFigFont{12}{14.4}{\rmdefault}{\mddefault}{\updefault}$\widetilde{A}_1^+$}}}
\put(3871,-3616){\makebox(0,0)[lb]{\smash{\SetFigFont{12}{14.4}{\rmdefault}{\mddefault}{\updefault}$\widetilde{A}_1^+$}}}
\put(4411,-4336){\makebox(0,0)[lb]{\smash{\SetFigFont{12}{14.4}{\rmdefault}{\mddefault}{\updefault}$\widetilde{A}_{t_+}$}}}
\put(2926,-4786){\makebox(0,0)[lb]{\smash{\SetFigFont{12}{14.4}{\rmdefault}{\mddefault}{\updefault}$\widetilde{A}_2^+$}}}
\put(2116,-4786){\makebox(0,0)[lb]{\smash{\SetFigFont{12}{14.4}{\rmdefault}{\mddefault}{\updefault}$\widetilde{A}_2^-$}}}
\put(6571,-4831){\makebox(0,0)[lb]{\smash{\SetFigFont{12}{14.4}{\rmdefault}{\mddefault}{\updefault}$p( \widetilde{A}_2^-)$}}}
\put(7381,-4831){\makebox(0,0)[lb]{\smash{\SetFigFont{12}{14.4}{\rmdefault}{\mddefault}{\updefault}$p(\widetilde{A}_2^+)$}}}
\put(8686,-4831){\makebox(0,0)[lb]{\smash{\SetFigFont{12}{14.4}{\rmdefault}{\mddefault}{\updefault}$p(\widetilde{A}_2^-)$}}}
\put(2881,-4336){\makebox(0,0)[lb]{\smash{\SetFigFont{12}{14.4}{\rmdefault}{\mddefault}{\updefault}$\widetilde{\R C}_{t_-}$}}}
\end{picture}
}}$$
It follows that
$$\begin{array}{rcl}
sp(\widetilde{\R C}_{t_+}) - sp(\widetilde{\R C}_{t_-}) &=& < \pi^* w (P_X) , 
p (\widetilde{\R C}_{t_+}) -
p (\widetilde{\R C}_{t_-})> \\
&=& < w (P_X) , \pi_* (p (\widetilde{A}_2^+) + p (\widetilde{A}_2^-)) > \quad (t_+, t_- \to t_1)\\
&=& 0,
\end{array}$$
where $w (P_X)$ is the characteristic class of the covering $P_X \to R_X$. $\square$

\begin{rem}
It would be interesting to extend Theorem \ref{maintheorem} to the case $\sum_{i=1}^n r_i = 0$.
The non-vanishing of this sum has been used to define the spinor state of real balanced 
rational
curves of class $d$ when $\R X$ is not orientable and $k_d = \frac{1}{2} c_1 (X) d$ is even. It
has also been used to make sure that all the real rational curves under
consideration have non-empty real parts. In the case of real curves with empty real parts, 
I cannot define a spinor state yet, and I 
cannot prevent the presence of critical points as the one described in Proposition \ref{propcrit}.
Note that this problem did not appear in dimension $2$, see \cite{Wels2}.
\end{rem}

\section{A study of the polynomial $\chi^{d, \s}$}
\label{sectfurth}

As was done in \cite{Wels2}, it is possible to understand how the invariant $\chi^{d, \s}_r$ 
depends on the $n$-tuple $r$ in terms of a new invariant. The aim of this last paragraph
is to explain this phenomenum in the case of $(X , c_X) = (\C P^3 , conj)$.

\subsection{Invariants of the blown up projective space}

\subsubsection{Statement of the results}
\label{subsectinvblownup}

Let $(Y , c_Y)$ be the real algebraic $3$-manifold obtained after blowing up a real point $x_0$
of $(\C P^3 , conj)$. Denote by $Exc \subset Y$ the exceptional divisor of the blow up,
$l \subset Exc \cong \C P^2$ the class of a line and $f \subset Y$ the strict transform
of a line in $\C P^3$ passing through $x_0$. The cone $NE (Y)$ of effective curves of $Y$
is of dimension $2$, closed and generated by $l$ and $f$. Note that $Y$ is not convex. Indeed,
if $u : \C P^1 \to Y$ has an image contained in $Exc$ which is not a line, then
$H^1 (\C P^1 ; u^* TY) \neq 0$. Note also that if $A \subset Y$ is an irreducible curve not
contained in $Exc$, then $A$ is equivalent to $af + bl$ with $0 \leq b \leq a$, since
$0 \leq A.Exc = a-b$.

Let $0 \leq k \leq d$ and $d_Y = d (f+l) - kl$, so that $c_1 (Y).d_Y = 4d - 2k$. Denote by
$k_{d_Y} = 2d-k$. Let $\underline{y} \subset Y^{k_{d_Y}}$ be a generic real configuration of
$k_{d_Y}$ distinct points of $Y$, and $r$ be the number of real points in this configuration
which is assumed to be non-zero. There are then only finitely many
connected rational curves in $Y$ passing through $\underline{y}$ in the homology class $d_Y$. 
Moreover, these curves are all irreducible and balanced and their
number does not depend on the generic choice of $\underline{y}$, it is equal
to the genus $0$ Gromov-Witten invariant $N_{d_Y} = GW (Y,d_Y, pt, \dots , pt)$ (see \cite{Gath}).
Denote by ${\cal R}_{d_Y} (\underline{y})$
the subset of these curves which are real and by  ${R}_{d_Y} (\underline{y})$ its
cardinality. Equip $\R Y = \R P^3 \# \overline{\R P}^3$ with a spin structure $\s_Y$. This spin 
structure allows us to define a spinor state for each real curve 
$ A \subset {\cal R}_{d_Y} (\underline{y})$, as in \S \ref{subsectspinorient}. Put then :
$$\chi^{d_Y,\s_Y}_r (\underline{y}) = \sum_{A \in {\cal R}_{d_Y} (\underline{y})} sp(A) \in \Z.$$
\begin{theo}
\label{theoremblownup}
Let $(Y , c_Y)$ be the blown up of $(\C P^3 , conj)$ at a real point $x_0$ and $\s_Y$ be a spin 
structure on $\R Y$. Let $0 \leq k \leq d$, $d_Y = d (f+l) - kl$ and $k_{d_Y} = 2d-k$.
Let $\underline{y}$ be a generic real 
configuration of $k_{d_Y}$ distinct points of $Y$ such that $r = \# (\underline{y} \cap \R Y)$
does not vanish. Then the integer 
$\chi^{d_Y,\s_Y}_r (\underline{y})$ does not depend on the choice of $\underline{y}$.
\end{theo}
This integer $\chi^{d_Y,\s_Y}_r (\underline{y})$ is denoted by $\chi^{d_Y,\s_Y}_r$. We deduce
from this Theorem \ref{theoremblownup} the following lower bounds in real enumerative
geometry.

\begin{cor}
\label{corlowerblown}
Under the assumptions of Theorem \ref{theoremblownup}, denote by ${R}_{d_Y} (\underline{y})$
the number of real irreducible rational curves passing through $\underline{y}$ in the
class $d_Y$ and by $N_{d_Y}$ the associated Gromov-Witten invariant. Then,
$|\chi^{d_Y,\s_Y}_r| \leq R_{d_Y} (\underline{y}) \leq N_{d_Y}. \quad \square$
\end{cor}

\subsubsection{Proof of Theorem \ref{theoremblownup}}

\begin{lemma}
\label{lemmanotconvex}
Let $Y$ be the blown up of $\C P^3$ at a point and $Exc \subset Y$ be the associated 
exceptional divisor. Then every morphism $u : \C P^1 \to Y$ such that $H^1 (\C P^1 ; u^* TY) \neq 0$
has an image contained in $Exc$ which is a curve of degree greater than one.
\end{lemma}

{\bf Proof :}

The manifold $Y$ has a submersion $p : Y \to \C P^2 = Exc$ whose fibers are rational curves. In fact,
$Y$ is isomorphic to the projectivization $P({\cal O}_{\C P^2} \oplus {\cal O}_{\C P^2} (1))$.
Denote by $H = {\cal O}_{\C P^2} \oplus {\cal O}_{\C P^2} (1)$, $L \subset p^*H$ the
tautological subline bundle over $Y$ and $M = p^*H/L$. Denote by $F$ the kernel of the morphism
$dp : TY \to T \C P^2$, one has the isomorphism $F \cong Hom (L , M) \cong L^* \otimes M$.
Let $u : \C P^1 \to Y$, from the short exact sequence 
$$0 \to {\cal O}_{\C P^1} (u^* F) \to {\cal O}_{\C P^1} (u^* TY) \to {\cal O}_{\C P^1} (u^* T\C P^2)
\to 0$$ we deduce the long exact sequence
$$\dots \to H^1 (\C P^1 ; u^* (L^* \otimes M)) \to H^1 (\C P^1 ; u^* TY) \to
H^1 (\C P^1 ; (p \circ u)^* (T\C P^2) \to \dots$$
Since $\C P^2$ is convex, $H^1 (\C P^1 ; (p \circ u)^* (T\C P^2) = 0$. Consider then the exact 
sequence  $0 \to u^* L \to (p \circ u)^* H \to u^* M \to 0$. Denote by $b$ the degree of
$p \circ u : \C P^1 \to \C P^2$, one has $(p \circ u)^* H \cong {\cal O}_{\C P^1} \oplus
{\cal O}_{\C P^1} (b)$. We deduce the following alternative. Either $u^* L$ is the unique
subline bundle of $(p \circ u)^* H$ of degree $b$, which implies $\deg ( u^* (L^* \otimes M))
= -b$ and $H^1 (\C P^1 ; u^* TY) \neq 0$ as soon as $b \geq 2$. In this case, $Im (u)$ is contained
in $Exc$. Or $u^* L$ is not the unique
subline bundle of $(p \circ u)^* H$ of degree $b$, then $\deg ( u^* (L)) \leq 0$ and
$\deg ( u^* (M)) \geq b$. This implies that $H^1 (\C P^1 ; u^*(L^* \otimes M) ) = 0$ and
thus $H^1 (\C P^1 ; u^* TY) = 0$. $\square$\\

{\bf Proof of Theorem \ref{theoremblownup} :}

Let $\tau \in {\cal S}_{k_{d_Y}}$ having $r$ fixed points in $\{ 1 , \dots , 
k_{d_Y} \}$ and such that $\tau^2 = id$. Let $(\overline{\cal M}^{d_Y}_{0 , k_{d_Y}} (Y),
c_{\overline{\cal M} , \tau})$ be the space of genus $0$ stable maps of $Y$ having $k_{d_Y}$
marked points and realizing the homology class $d_Y$. This is a real algebraic subvariety
of the variety $(\overline{\cal M}^{b}_{0 , k_{d_Y}} (\C P^N),
c_{\overline{\cal M} , \tau})$ associated to a real embedding $\pi : (Y, c_Y) \to
(\C P^N , conj)$ such that $\pi_* d_Y = b$, see \cite{FP}. Denote by 
$\overline{\cal M}^{d_Y}_{0 , k_{d_Y}} (Y)^\#$ the complement in 
$\overline{\cal M}^{d_Y}_{0 , k_{d_Y}} (Y)$ of the locus of curves having a multiple component or
a component included in $Exc$. This variety is smooth of dimension $c_1 (Y) d_Y + k_{d_Y} =
3k_{d_Y}$, but not compact. However, the evaluation morphism 
$\overline{\cal M}^{d_Y}_{0 , k_{d_Y}} (Y)^\# \to Y^{k_{d_Y}}$ is proper over a generic path
of $Y^{k_{d_Y}}$. Theorem \ref{theoremblownup} is then proven in the same way as 
Theorem \ref{maintheorem}. $\square$

\subsection{Relations between the coefficients of the polynomial $\chi^{d, \s} (T)$}
\label{subsectfurth}

Denote once more by $(Y , c_Y)$ the real algebraic $3$-manifold obtained after blowing up a real 
point $x_0$ of $(\C P^3 , conj)$. Assume that the riemannian metric $g$ of $\R P^3$ fixed
in \S \ref{subsectspinstruct} is flat in a neighborhood of $x_0$. Denote by $\overline{g}$
the associated metric of $\overline{\R P}^3$ and by $\overline{R}_X$ the $SO_3 (\R)$-principal
bundle of direct orthonormal frames of $\overline{\R P}^3$. Let $\overline{P}_X \to \overline{R}_X$
be a $Spin_3$-principal bundle defining a spin structure $\overline{\s}$ on $\overline{\R P}^3$.
Making a surgery in a neighborhood of $x_0$ which is small enough, the manifold $\R Y =
\R P^3 \# \overline{\R P}^3$ comes equipped with a riemannian metric $g \# \overline{g}$.
Denote by $R_Y$ the associated $SO_3 (\R)$-principal
bundle of direct orthonormal frames. The bundles $P_X$ and $\overline{P}_X$ glue together
to form a $Spin_3$-principal bundle $P_Y$ over $R_Y$ which defines a spin structure on $\R Y$
denoted by $\s \# \overline{\s}$.
\begin{theo}
\label{theorel}
Let $(Y , c_Y)$ be the blown up of $(\C P^3 , conj)$ at a real point $x_0$, $\s$ be a spin 
structure on $\R P^3$ and $\s \# \overline{\s}$ be the associated  spin 
structure on $\R Y$. Let $d \in \N^*$, $d_Y = d (f+l) - 2l \in H_2 (Y ; \Z)$ and $k_{d} = 2d$.
Then, for every integer $r$ between $2$ and $k_d -2$, we have 
$\chi^{d,\s}_{r+2} = \chi^{d,\s}_{r} - 2\chi^{d_Y,\s \# \overline{\s}}_{r}.$
\end{theo}
(The classes $f$ and $l$ in $ H_2 (Y ; \Z)$ have been defined in \S \ref{subsectinvblownup}.)

\begin{lemma}
\label{lemmablowup}
Let $u : \C P^1 \to \C P^3$ be a morphism passing through $x_0 \in \C P^3$, which is an immersion 
in the neighborhood of $u^{-1} (x_0)$. Denote by $N_u$ the normal bundle of $u$ in $\C P^3$ and
by $z = u^{-1} (x_0) \subset \C P^1$. Let $Y$ be the blown up of $\C P^3$ at $x_0$,
$\tilde{u} : \C P^1 \to Y$ be the strict transfrom of $u$ and $N_{\tilde{u}}$ be its normal bundle
in $Y$. Then one has the isomorphism :
$$N_{\tilde{u}} \cong N_u \otimes {\cal O}_{\C P^1} (-z). \quad \square$$
\end{lemma}
It follows in particular from this lemma that if $u$ is balanced, then so is $\tilde{u}$.\\

{\bf Proof of Theorem \ref{theorel} :}

The begining of this proof is the same as the one of Theorem $3.2$ of \cite{Wels2}. Let 
$\underline{x}^+ = (x_1 , \dots , x_{k_d})$ be a generic real configuration of $k_d$ distinct
points of $\C P^3 \setminus \{ x_0 \}$, $r+2$ of which are real. Let two such real points,
say $x_{k_d - 1}$ and $x_{k_d}$, converge to $x_0$ along a generic real tangency $\tau_0 \in
P(T_{x_0} \R P^3)$. Denote by $\underline{x}_{\infty}$ the configuration of points
$(x_1 , \dots , x_{k_d -2} , x_0)$ in the limit. The rational curves in the set 
${\cal R}_d (\underline{x}^+)$
converge to immersed irreducible real rational curves passing through $\underline{x}_{\infty}$ and
either having the tangency $\tau_0$ at $x_0$, or having a real ordinary double point at $x_0$
which is the local intersection of two real branches. The latter are the limit of exactly 
two curves of ${\cal R}_d (\underline{x}^+)$, denote by ${\cal L}im^+ (\underline{x}_{\infty})$
their set. Similarly, $\underline{x}_{\infty}$ is the limit of a generic real configuration
$\underline{x}^-$ of $k_d$ distinct
points of $\C P^3 \setminus \{ x_0 \}$, $r$ of which being real, when two complex conjugated points
of $\underline{x}^-$ converge to $x_0$ along a generic real tangency $\tau_0 \in
P(T_{x_0} \R P^3)$. The rational curves in the set ${\cal R}_d (\underline{x}^-)$
converge to immersed irreducible real rational curves passing through $\underline{x}_{\infty}$ and
either having the tangency $\tau_0$ at $x_0$, or having a real ordinary double point at $x_0$
which is the local intersection of two complex conjugated branches. The latter are the limit of 
exactly 
two curves of ${\cal R}_d (\underline{x}^-)$, denote by ${\cal L}im^- (\underline{x}_{\infty})$
their set. It follows easily that
$$\chi^{d , \s}_{r+2} - \chi^{d , \s}_{r} = 2\Big( \sum_{A \in {\cal L}im^+ (\underline{x}_{\infty})}
sp(A) - \sum_{A \in {\cal L}im^- (\underline{x}_{\infty})} sp(A) \Big).$$
For every curve $A \in {\cal L}im (\underline{x}_{\infty}) = {\cal L}im^+ (\underline{x}_{\infty})
\cup {\cal L}im^- (\underline{x}_{\infty})$, denote by $A_Y$ its strict transform in $Y$. The
curves $A_Y$, $A \in {\cal L}im (\underline{x}_{\infty}) $, are exactly the elements of the set
${\cal R}_{d_Y} (\underline{y})$, for $\underline{y} = (x_1 , \dots , x_{k_d -2}) \in Y^{k_{d_Y}}$.
When $A \in {\cal L}im^- (\underline{x}_{\infty})$, it follows from the construction of the spin
structure $\s \# \overline{\s}$ that $sp(A_Y) = sp(A)$. We have thus to prove that when 
$A \in {\cal L}im^+ (\underline{x}_{\infty})$, $sp(A_Y) = -sp(A)$.

Let $(e_1 (p) , e_2 (p), e_3 (p))_{p \in \R A_Y}$ be a loop in the bundle $R_Y$ of direct orthonormal
frames of $\R Y$ constructed in \S \ref{subsectspinorient}, so that $sp (A_Y)$ is the obstruction
to lift this loop as a loop in $P_Y$. Fix an orientation of $\R A_Y$, which induces also an 
orientation of the normal bundle of $\R A_Y$ in $\R Y$, and let 
$\R A_Y \cap \R Exc = \{ p_0 , p_1 \}$. We construct then the loop $(\tilde{e}_1 (p) , 
\tilde{e}_2 (p), \tilde{e}_3 (p))$ of $R_Y$ as the concatenation of the four following paths.
The path $(e_1 (p) , e_2 (p), e_3 (p))_{p \in [p_0 , p_1]}$, then a path completely included
in the fiber of $R_Y$ over $p_1$. This path is obtained from $(e_1 (p_1) , e_2 (p_1), e_3 (p_1))$
by having this frame turning of half a twist in the
positive direction around the axis generated by $e_1 (p_1)$. The end point of this path is
thus the frame $(e_1 (p_1) , -e_2 (p_1), -e_3 (p_1))$. Then the path 
$(e_1 (p) , -e_2 (p), -e_3 (p))_{p \in [p_1 , p_0]}$. Finally, a  path completely included
in the fiber of $R_Y$ over $p_0$. This path is obtained from $(e_1 (p_0) , -e_2 (p_0), -e_3 (p_0))$
by having this frame turning of half a twist in the
negative direction around the axis generated by $e_1 (p_0)$. The end point of this path is
thus the frame $(e_1 (p_0) , e_2 (p_0), e_3 (p_0))$, so that we have indeed constructed a loop
in $R_Y$ which is homotopic to the initial loop $(e_1 (p) , e_2 (p), e_3 (p))_{p \in \R A_Y}$.
Remember that $\R Y$ is obtained as the connected sum of $\R P^3$ and $\overline{\R P}^3$ in the
neighborhood of $x_0$. We may choose the radius of the ball used to perform the connected sum
as small as we want, and have this radius converging to zero. In this process, the curve $\R A_Y$
degenerates to the union of the curve $\R A \subset \R P^3$ and two lines $\R D_1, \R D_2$ of
$\overline{\R P}^3$ passing through $x_0$. The loop $(\tilde{e}_1 (p) , 
\tilde{e}_2 (p), \tilde{e}_3 (p))$ of $R_Y$ degenerates to the union of a loop
$(\tilde{e}_1  ,  \tilde{e}_2 , \tilde{e}_3)|_{\R A}$ of $R_X$ and two loops 
$(\tilde{e}_1  ,  \tilde{e}_2 , \tilde{e}_3)|_{\R D_1}$, 
$(\tilde{e}_1  ,  \tilde{e}_2 , \tilde{e}_3)|_{\R D_2}$ of $\overline{R}_X$. By construction,
the obstruction to lift $(\tilde{e}_1  ,  \tilde{e}_2 , \tilde{e}_3)|_{\R A}$ to $P_X$ is
exactly $sp(A)$. Similarly, the homotopy class of the loop 
$(\tilde{e}_1  ,  \tilde{e}_2 , \tilde{e}_3)|_{\R D_1}$ differs from the one of
$(\tilde{e}_1  ,  \tilde{e}_2 , \tilde{e}_3)|_{\R D_2}$ by a non-trivial loop in a fiber of
$\overline{P}_X$. We deduce that $sp(\R D_1) = - sp(\R D_2)$ independantly of the choice of 
a spin structure on $\overline{\R P}^3$. It follows that
$sp(A_Y) = sp(A) sp(\R D_1) sp(\R D_2) = -sp(A) sp(\R D_1)^2 = -sp(A). \quad \square$

\addcontentsline{toc}{part}{\hspace*{\indentation}Bibliography}

 \nocite{*}  

\bibliography{convexe}

\begin{thebibliography}{10}

\bibitem{ABS}
M.~F. Atiyah, R.~Bott, and A.~Shapiro.
\newblock Clifford modules.
\newblock {\em Topology}, 3(suppl. 1):3--38, 1964.

\bibitem{FP}
W.~Fulton and R.~Pandharipande.
\newblock Notes on stable maps and quantum cohomology.
\newblock In {\em Algebraic geometry---Santa Cruz 1995}, volume~62 of {\em
  Proc. Sympos. Pure Math.}, pages 45--96. Amer. Math. Soc., Providence, RI,
  1997.

\bibitem{Gath}
A.~Gathmann.
\newblock Gromov-{W}itten invariants of blow-ups.
\newblock {\em J. Algebraic Geom.}, 10(3):399--432, 2001.

\bibitem{Gro}
M.~Gromov.
\newblock Pseudoholomorphic curves in symplectic manifolds.
\newblock {\em Invent. Math.}, 82(2):307--347, 1985.

\bibitem{IKS}
I.~Itenberg, V.~Kharlamov, and E.~Shustin.
\newblock Welschinger invariant and enumeration of real rational curves.
\newblock {\em Internat. Math. Res. Notices}, 49:2639--2653, 2003.

\bibitem{KirT}
R.~C. Kirby and L.~R. Taylor.
\newblock {${\rm Pin}$} structures on low-dimensional manifolds.
\newblock In {\em Geometry of low-dimensional manifolds, 2 (Durham, 1989)},
  volume 151 of {\em London Math. Soc. Lecture Note Ser.}, pages 177--242.
  Cambridge Univ. Press, Cambridge, 1990.

\bibitem{Ko}
J.~Koll{\'a}r.
\newblock {\em Rational curves on algebraic varieties}, volume~32 of {\em
  Ergebnisse der Mathematik und ihrer Grenzgebiete. 3. Folge. A Series of
  Modern Surveys in Mathematics [Results in Mathematics and Related Areas. 3rd
  Series. A Series of Modern Surveys in Mathematics]}.
\newblock Springer-Verlag, Berlin, 1996.

\bibitem{Kont}
M.~Kontsevich and Y.~Manin.
\newblock Gromov-{W}itten classes, quantum cohomology, and enumerative
  geometry.
\newblock {\em Comm. Math. Phys.}, 164(3):525--562, 1994.

\bibitem{Kwon}
S.~Kwon.
\newblock Real aspects of the moduli space of stable maps of genus zero curves.
\newblock {\em Preprint math.AG/0305128}, 2003.

\bibitem{Miln}
J.~W. Milnor and J.~D. Stasheff.
\newblock {\em Characteristic classes}.
\newblock Princeton University Press, Princeton, N. J., 1974.
\newblock Annals of Mathematics Studies, No. 76.

\bibitem{OSS}
C.~Okonek, M.~Schneider, and H.~Spindler.
\newblock {\em Vector bundles on complex projective spaces}, volume~3 of {\em
  Progress in Mathematics}.
\newblock Birkh\"auser Boston, Mass., 1980.

\bibitem{RT}
Y.~Ruan and G.~Tian.
\newblock A mathematical theory of quantum cohomology.
\newblock {\em J. Differential Geom.}, 42(2):259--367, 1995.

\bibitem{Sot}
F.~Sottile.
\newblock Enumerative real algebraic geometry.
\newblock {\em Algorithmic and Quantitative Aspects of Real Algbraic Geometry,
  S. Basu and L. Gonzalez-Vega, eds., DIMACS series 60, AMS}, 2003.
\newblock To appear.

\bibitem{Wels}
J.-Y. Welschinger.
\newblock Invariants of real rational symplectic 4-manifolds and lower bounds
  in real enumerative geometry.
\newblock {\em C. R. Acad. Sci. Paris S\'er. I Math.}, 336(4):341--344, 2003.

\bibitem{Wels2}
J.-Y. Welschinger.
\newblock Invariants of real symplectic 4-manifolds and lower bounds in real
  enumerative geometry.
\newblock {\em Pr\'epublication de l'Ecole Normale Sup\'erieure de Lyon}, 312,
  2003.

\bibitem{Wit}
E.~Witten.
\newblock Two-dimensional gravity and intersection theory on moduli space.
\newblock In {\em Surveys in differential geometry (Cambridge, MA, 1990)},
  pages 243--310. Lehigh Univ., Bethlehem, PA, 1991.

\end{thebibliography}
\bibliographystyle{abbrv}

\noindent Ecole Normale Sup\'erieure de Lyon\\
Unit\'e de Math\'ematiques Pures et Appliqu\'ees\\
UMR CNRS $5669$\\
$46$, all\'ee d'Italie\\
$69364$, Lyon Cedex $07$\\
(FRANCE)\\
e-mail : {\tt jwelschi@umpa.ens-lyon.fr}

\end{document}